\documentclass[11pt]{amsart}
\usepackage{amssymb, amsmath, amsthm, amsfonts}
\usepackage{epsfig}
\usepackage{latexsym,amssymb}
\numberwithin{equation}{section}
\newtheorem{teo}{Theorem}[section]
\newtheorem{cor}[teo]{Corollary}
\newtheorem{prop}[teo]{Proposition}
\newtheorem{lemma}[teo]{Lemma}

\newtheorem{remark}[teo]{Remark}
\newtheorem{example}[teo]{Example}
\newtheorem{definition}[teo]{Definition}

\theoremstyle{definition}

\theoremstyle{remark}

\font\maius=cmcsc10 scaled1200
\font\svfilt=msbm7

 \font\filt=msbm10
\def\neweq#1{\begin{equation}\label{#1}}
\def\endeq{\end{equation}}

\def\proof{\noindent{\it Proof.} }

\def\A{\mathcal A}
\def\C{\mathcal C}
\def\D{\delta J (f, g)}
\def\g{\gamma_n}

\def\dom{{\rm dom}}
\def\Ent{{\rm Ent}}
\def\e{\varepsilon}
\def\H{{\mathcal H}}
\def\K{\mathcal K}
\def\L{\mathcal L}
\def\ov{\overline}

\def \R {\hbox{\filt R}}
\def \res{\mathop{\hbox{\vrule height 7pt width .5pt depth 0pt\vrule height .5pt width 6pt depth 0pt}}\nolimits}
\def\sfe{S^{n-1}}
\def \sR {\hbox{\svfilt R}}

\begin{document}

\title[The area measure of log-concave functions and related inequalities]{The area measure of log-concave functions
\\ and related inequalities}
%\date{\today}
\author{Andrea Colesanti \ \& \ Ilaria Fragal\`a}
\address{Andrea Colesanti, Dipartimento di Matematica ``U. Dini'', Universit\`a degli Studi di Firenze, Viale Morgagni 67/A, 50134 Firenze (Italy)}
\address{Ilaria Fragal\`a, Dipartimento di Matematica, Politecnico di Milano, Piazza Leonardo Da Vinci 32, 20133 Milano (Italy)}

\baselineskip14pt

\maketitle
\begin{abstract} On the class of log-concave functions on
$\R^n$, endowed with a suitable algebraic structure, we study the
first variation of the total mass functional, which corresponds to
the volume of convex bodies when restricted to the subclass of
characteristic functions. We prove some integral representation
formulae for such first variation, which lead to define in a natural
way the notion of area measure for a log-concave function. In the
same framework, we obtain a functional counterpart of Minkowski
first inequality for convex bodies; as corollaries, we derive a
functional form of the isoperimetric inequality, and a family of
logarithmic-type Sobolev inequalities with respect to log-concave
probability measures. Finally, we propose a suitable functional
version of the classical Minkowski problem for convex bodies, and
prove some partial results towards its solution.
\end{abstract}

\medskip
\noindent {\small {\it 2010MSC$\,$}: 26B25 (primary), 26D10, 52A20.}

\smallskip
\noindent {\small {\it Keywords$\,$}: Log-concave functions, convex
bodies, area measure, isoperimetric inequality, log-Sobolev
ine\-qua\-lit\-y, Minkowski problem.}

\bigskip\bigskip \centerline{{\maius Plan of the paper}}

1. Introduction

2. Preliminaries

\ \ 2.1 Notation and background

\ \ 2.2 Functional setting

3. Differentiability of the total mass functional

\ \ 3.1 Existence of the first variation

\ \ 3.2 Computation of the first variation in some special cases

4. Integral representation of the first variation

5. The functional form of Minkowski first inequality

6. Isoperimetric and log-Sobolev inequalities for log-concave
functions

7. About the Minkowski problem

8. Appendix

\section{Introduction}

This article regards {\it log-concave} functions defined in $\R^n$,
{\it i.e.} functions of the form
$$
f\,:\,\R^n\,\rightarrow\,\R\,,\quad f=e^{-u}\,,
$$
where
$
u:\R^n\,\rightarrow\,\R\cup\{+\infty\}
$
is convex.

In the last decades the interest for log-concave functions has been
considerably increasing, strongly motivated by the analogy between
these objects and convex bodies (convex compact subsets of $\R^n$).

The first breakthrough in the discovery of parallel behaviours of
convex bodies and log-concave functions, was the {\it
Pr\'ekopa-Leindler inequality}, named after the two Hungarian
mathematicians who proved it in the seventies \cite{Le,Pr1,Pr2,
Pr3}. It states that, for any given functions $f, g,h \in L ^ 1 (\R
^n; \R _+)$ which satisfy, for some $t \in (0,1)$, the pointwise
inequality
$$h \big ((1- t ) x + t y \big ) \geq f(x) ^ {1 - t} g ( y ) ^ t
\qquad \forall x, y \in \R ^n\ ,
$$
it holds
\begin{equation}\label{intro1}
\int _{\sR ^n} h \geq \Big ( \int _{\sR ^n} f   \Big ) ^ {1- t} \Big
( \int _{\sR ^n} g   \Big ) ^ {t} \ .
\end{equation}
Moreover, it was proved by Dubuc in \cite{Dubuc} that  the equality
sign holds in (\ref{intro1}) if and only if the functions $f$ and
$g$ are log-concave and translates, meaning that $f (x) = g ( x-
x_0)$ for some $x_0 \in \R ^n$.

If $K$ and $L$ are measurable subsets of $\R^n$ such that also their Minkowski combination
$(1-t) K + tL$ is measurable, by
applying the Pr\'ekopa-Leindler inequality with  $f, g$ and $h$ equal respectively to the characteristic functions of $K$, $L$ and $(1-t) K + tL$, one obtains
$$
V((1-t)K+tL)\ge V(K)^{1-t}\,V(L)^t\,.
$$
This is an equivalent formulation of the classical {\it Brunn-Minkowski inequality}
%(see \cite{Gardner} or \cite{Schneider}):
\begin{equation}\label{intro2}
V((1-t)K+tL)^{1/n}\ge (1-t)V(K)^{1/n} +t V(L)^{1/n}\,,
\end{equation}
which holds with equality sign if and only if $K$ and $L$ belong to
the class $\K ^n$ of convex bodies in $\R^n$ and are homothetic,
namely they agree up to a translation and a dilation.

The geometric inequality (\ref{intro2}) is a cornerstone in Convex
Geometry: it has many important consequences, such as the
isoperimetric inequality for convex bodies, and the uniqueness issue
in the solution of the Minkowski problem (see the survey paper
\cite{Gardner} for an overview). On the other hand, in view of its
functional form, inequality (\ref{intro1}) is somehow more
``flexible'', and finds many applications in different fields, such
as convex geometry, probability, mass transportation; we refer the
reader to \cite{Ball,Barthe,Villani} for more information on
Pr\'ekopa-Leindler inequality, including proofs and bibliographical
references.

In the same way as (\ref{intro1}) paraphrases
(\ref{intro2}) into the realm of functions, recently analytic versions of
other geometric inequalities have been studied. In particular, we
mention the so-called Blaschke-Santal\'o inequality, involving the
product of the volume of a convex body and its polar: functional
versions of it have been achieved by Ball \cite{Ball}, Artstein,
Klartag and Milman \cite{Artstein-Klartag-Milman04}, and Fradelizi
and Meyer \cite{Fradelizi-Meyer}. Let us also emphasize that a
suitable notion of mean width for log-concave functions has been
introduced by Klartag and Milman in  \cite{Klartag-Milman05}, where
some related Urysohn-type inequality are also proved; a short time
ago, these topics have been further developed by Rotem in
\cite{Rotem1, Rotem2}.

\medskip
In the same spirit, the aim of this paper is to cast some more light
upon the geometry of log-concave functions, and to propose
functional counterparts of some classical quantities and
inequalities in Convex Geometry, that we briefly remind below (for
more details, we refer to \cite{Schneider}).

Going back to the Brunn-Minkowsi inequality, let us recall that
it admits a sort of ``differential version'', the so-called
{\it Minkowski first inequality}, which reads
\begin{equation}\label{intro3}
V_1(K,L):= \frac{1}{n}\,\lim_{t\to0^+}
\frac{V(K+tL)-V(L)}{t} \ge V(K)^{\frac{n-1}{n}}V(L)^{\frac{1}{n}} \qquad \forall K, L \in \K ^n\ .
\end{equation}

Inequality (\ref{intro3}) can be easily obtained from
(\ref{intro2}), and it is in fact equivalent to it.
%(see for instance
%\cite[Section 7]{Gardner})
Notice that, when $L$ is the unit ball, $V_1(K,L)$ is
just the perimeter of $K$, and (\ref{intro3})
becomes the isoperimetric inequality in the class of convex bodies.

The term $V_1(K,L)$, which is one
of the mixed volumes of $K$ and $L$,
%(see \cite[Chapter 5]{Schneider})
admits a very simple and elegant integral representation:
\begin{equation}\label{intro5}
V_1(K,L)=\frac{1}{n}\int_{{\sfe}} h_L\,d\sigma_K\,,
\end{equation}
where $h_L$ is the {\it support function} of $L$,  and $\sigma _K$ is the {\it area measure} of $K$.
%The area measure $\sigma_K$ arises naturally
%as one of the coefficient measures in the local Steiner formula for $K$ (see
%\cite[Chapter 4]{Schneider}), but it can be equivalently identified by
%(\ref{intro5}) as the unique measure such that this equality holds true for
%every $L$.
In view of (\ref{intro5}), the measure $\sigma_K$ is usually
interpreted as the first variation of volume with respect to the
Minkowski addition. The classical {\it Minkowski problem} consists
in retrieving $K$ from its surface area measure, and it is
well-known that it admits a unique solution up to translations. More
precisely, given any measure $\eta$ on the unit sphere $S ^ {n-1}$
which satisfies the compatibility conditions of having null
barycenter and being not concentrated on an equator, there exists a
convex body, unique up to translations, such that $\eta = \sigma
_K$.

\medskip
Our main goals are to provide a functional version of Minkowski
first inequality (\ref{intro3}), of the representation formula
(\ref{intro5}), and of the Minkowski problem. In this perspective, a
crucial issue is to identify a good notion of ``area measure'' for a
log-concave function. To that aim, we pursue a quite natural idea,
namely we replace the volume of a convex body by the integral of a
log-concave function: we set
$$J (f) = \int_{\sR^n} f \, dx\ ,$$
and we compute the  first variation of such integral functional with respect to suitable perturbations.

Actually, log-concave functions can be equipped with two internal operations:
a sum and a multiplication by positive reals, that will be denoted respectively by
$\oplus$ and $\cdot$, and can be
characterized as follows (see Section 2 for a more rigorous presentation). If $f=e^{-u}$ and $g=e^{-v}$ are
log-concave functions and $\alpha, \beta>0$, then
\begin{equation}
\label{intro00}
\alpha\cdot f\oplus\beta\cdot g:=e^{-w}\,,
\quad
\mbox{where $w^*= \alpha u^*+\beta v^*$\ .}
\end{equation}
Here ${}^*$ denotes as usual the Fenchel conjugate of convex
functions. In other words, if we write a generic
log-concave function as $e^{-u}$, the operations introduced in (\ref{intro00}) are
linear with respect to $u^*$.
In particular, since the Fenchel conjugate of the indicatrix of a convex body is precisely its support function,  one has
$$
\alpha\cdot\chi_K\oplus\beta\cdot\chi_L=\chi_{\alpha K+\beta L}\ .
$$
Therefore, definition (\ref{intro00}) can be seen a natural extension to the class log-concave functions of  the Minkowski structure on convex bodies.

In this framework, for a pair of log-concave functions $f$ and $g$,
we study the quantity
\begin{equation}\label{intro6}
\delta J(f,g):=\lim_{t\to0^+}\frac{J(f\oplus t\cdot g)-J(f)}{t}\,.
\end{equation}

Let us point out that, red within this formalism, the above quoted works \cite{Klartag-Milman05,Rotem1,Rotem2} are concerned precisely with the limit in (\ref{intro6}), in the special case when $f$ is equal to $\gamma _n$, the density of the Gaussian
measure in $\R^n$. In fact, to some extent, $\g$ plays the role of the unit ball in the class of
log-concave functions. Thus, according to \cite{Klartag-Milman05},
the {\it mean width} of a
log-concave function $g$  is given by $\delta J(\g,g)$, by analogy with
the mean width of a convex body $K$ which is given by $V_1(B,K)$.
We also mention the paper \cite{Klartag} by Klartag (see in particular \S 3), where a limit similar to 
(\ref{intro6}) is considered, in the class of $s$--concave functions endowed with the 
appropriate algebraic operations, in order to derive several functional inequalities.

%In \cite{Klartag-Milman05} it is proved a
%Urysohn type inequality for $M(g)$, while in \cite{Rotem2} the
%author establishes a representation formula, which turns out to be a
% special case of Theorem \ref{introA'} of the present paper.

\medskip
When $f$ and $g$ are arbitrary log-concave functions, the limit in
(\ref{intro6}) exists under the fairly weak condition $J (f)>0$. In
Section 3.1 we give a rigorous proof of this fact, already pointed
out in \cite{Klartag-Milman05}, and we show that the condition $J
(f)>0$ is not necessary in the one dimensional case. Moreover we
give simple examples which reveal that $\delta J(f,g)$  may become
negative or $+\infty$ (indeed, whereas $V(K+tL)$ is a polynomial in
$t$ for every $K$ and $L$ in $\K ^n$, this is no longer true in
general for $J(f\oplus t\cdot g)$). Then in Section 3.2 we compute
$\delta J (f, g)$ in some special cases: the case when $f = g$,
which brings into play the {\it entropy} of $f$:
$$
\Ent(f) = \int_{\sR^n} f \log f \, dx - J (f) \log J (f)\ ,
$$
and the case when the logarithms of $f$ and $g$ are powers of
support functions of convex bodies, which allows to recover an
integral representation formula for the derivative of $p$-mixed
volume due to Lutwak \cite{Lu}.

\medskip
To go farther than these special cases, in Section 4 we come to the
problem at the core of the paper, namely the problem of giving some
general integral representation formula for $\delta J(f,g)$. We are
able to achieve such a representation in two distinct settings: when
the finiteness domains of $u = -\log f$ and $v= -\log g$ are the
whole space $\R^n$, and when such domains are smooth strictly convex
bodies. In both cases we have to assume further properties on $u$
and $v$, concerning regularity, growth at the boundary of their
domain, and strict convexity. To be more precise, our integral
representation formulae are settled in the classes $\A'$, $\A''$ of
log-concave functions $f= e ^ {-u}$ such that $u$ belongs
respectively  to
$$
\begin{array}{ll}
&
 \L ' := \Big \{ u \in \L \ :\ \dom (u)= \R^n
,\quad
 u \in \C ^ 2 _+ (\R ^n) \ ,\quad  \lim \limits_{\|x\| \to + \infty} \frac{u (x)} {\|x\|} = +
\infty \Big \}\,,
\\  \noalign{\medskip}
& \L''  := \Big \{ u \in \L\ :\  \dom (u) = K \in \K ^n \cap \C ^2 _+, \quad u
\in \C^2_+ ({\rm int} (K) ) \cap \C ^0 (K)\ ,\quad \lim \limits_{x
\to
\partial K} \| \nabla u (x) \| = + \infty \Big \}
\ .
\end{array}
$$

Here the notation $\C^2_+$, used for functions and sets, has the
following standard meaning: when it is referred to a function $u$,
it means that $u\in\C^2$ and the Hessian matrix of $u$ is positive
definite at each point; when it is referred to a convex body $K$, it
means that $\partial K\in\C^2$ and the Gauss curvature is everywhere
strictly positive.

After proving that $\A'$ and $\A''$ are both closed with respect to
the operations $\oplus$ and $\cdot$ (see Lemma \ref{chiusura}), we
state our main results, which are valid under the assumption that
the perturbation $g$ is ``controlled'' by the perturbed function $f$
(see Definition \ref{defH} for the precise statement of this
assumption, which is not necessary in the one dimensional case).  In
Theorem \ref{A'} we prove that, when $f,g \in \A'$, $\D$ is finite
and is given by
\begin{equation}\label{introtesi'}
\D = \int _{\sR^n}v^*(\nabla u(x))f(x)\,dx\ .
\end{equation}
 In
Theorem \ref{A''} we prove that, when $f,g \in \A''$, $\D$ is finite
and is given by
\begin{equation}\label{introtesi''} \D=\int _{K}v^*(\nabla
u(x))f(x)\,dx+ \int_{\partial K} h_L(\nu_K(x))f(x)\,d\H^{n-1}\,,
\end{equation}
where $K=\dom(u)$, $\nu_K$ is the unit outer normal to $\partial K$,
$L=\dom(v)$, and $h _L$ is the support function of $L$. The proof of
these results is quite delicate and requires a careful analysis, see
Section 4.

If we perform the change of variable $\nabla u(x)=y$ in
(\ref{introtesi'}), it becomes
\begin{equation}\label{introduale1}
\D=\int_{\sR^n}v^* \, d \mu (f) \ , \qquad   d \mu (f) := f (y)
 e^{-\langle y, \nabla u^*(y)\rangle +u^*(y)} \, \det (\nabla ^
2u^*(y)) \, dy\ .
\end{equation}
Comparing (\ref{introduale1}) with (\ref{intro5}), we are lead to
identify the measure $\mu(f)$ as the {\it area measure} of a
function $f$ in the class $A'$. (Under this point of view,  $v^*$
plays the role of support function of $g$, as in
\cite{Klartag-Milman05}; this interpretation is quite natural if we
remind that the algebraic structure we put on log-concave functions
$e^{-u}$ is linear with respect to $u^*$, in the same way as the
Minkowski structure on $\K ^n$ is linear with respect to support
functions). Similarly, with the changes of variable $\nabla u(x)=y$
and $\nabla \nu _K (y) = \xi$, (\ref{introtesi''}) becomes
\begin{equation}\label{introduale2}
\D=\int_{\sR^n}v^* \, d \mu (f)  + \int_{S^{n-1}} h _L \, d \sigma
(f) \ , \quad d \mu (f)  \hbox{ as above,} \quad d\sigma(f)
:=f(\nu_K^{-1}(\xi))\,d\sigma_K(\xi)  \ .
\end{equation}
Hence, within the class $\A''$, the notion of area measure of $f$ is
provided by the pair  $(\mu(f),\sigma(f))$ (notice that the former
is a measure on $\R^n$, the latter on $S ^ {n-1}$).

\medskip
Having the above representation formulae at our disposal, we then
turn attention to functional inequalities involving $\D$. 
Our approach is similar 
to the one used by Klartag in \cite{Klartag} for the class of $s$-concave functions.
In Section 5, we prove the following functional form of Minkowski first
inequality (\ref{intro3}) (see Theorem \ref{teomink}):
\begin{equation}\label{intromink1}
\delta J (f, g) \geq J (f) \big [\log  {J (g)}
 +n \big ] + \Ent (f) \ ,
\end{equation}
with equality sign if and only if there exists $x_0\in \R ^n$ such
that $g(x) = f (x-x_0)$ $\forall x \in \R^n$. Loosely speaking,
(\ref{intro3}) can be proved taking the right derivative at $t=0$ of
both sides of the Brunn-Minkowski inequality (\ref{intro2}), and
inequality (\ref{intromink1}) is obtained by adapting this idea to
the Pr\'ekopa-Leindler inequality, and using Dubuc's
characterization of the equality case.

In Section 6 we show that, by combining the abstract inequality
(\ref{intromink1}) with the above representation formulae for $\D$,
further functional inequalities come out.

Firstly, we define the {\it perimeter} of a function $f \in \A'$ in
the natural way, that is as $P (f):= \delta J (f, \gamma _n)$, and
we show that, under suitable assumptions, the following functional
version of the isoperimetric inequality holds (see Proposition
\ref{repper}):
$$P (f) = \frac{1}{2} \int_{\sR ^n} \frac{ \|\nabla f\| ^ 2} { f} \, dx + (\log
c_n ) \, J (f) \geq n J (f) + \Ent (f) \ ,
$$
with equality sign if and only if there exists $x_0 \in \R ^n$ such
that $f (x) = \g  (x- x_0)$ $\forall x \in \R ^n$.

Then we derive a family of inequalities of logarithmic Sobolev type
for probability measures $\nu$ with log-concave densities: under
suitable assumptions on $\nu$, $a$ and $h$, we obtain (see
Proposition \ref{logsobolev})
\begin{equation}\label{introLS'}
\int_{\sR ^n} a(h) \log a (h)
 d \nu   -  \Big ( \int_{\sR ^n} a(h) \, d \nu  \Big ) \log   \Big ( \int_{\sR ^n} a(h) \, d \nu  \Big )
 \leq  \frac{1}{c}
\int_{\sR ^n} \frac{(a' (h)) ^ 2} {a(h)} \|\nabla h \| ^ 2 \, d \nu
  \ .
 \end{equation}
In particular, by choosing  $\nu = \g \, dx$ and $a(h) = h ^2$, we
recover Gross' logarithmic Sobolev inequality for the Gaussian
measure.  We point out that our approach allows much more general
choices of $\nu$ and $a$; on the other hand, as a drawback, the
validity of (\ref{introLS'}) is obtained under some further
restrictions on $h$.

\medskip

Finally, in Section 7 we move few steps towards the solution of the
Minkowski problem for log-concave functions. As a natural extension
of the Minkowski problem for convex bodies, such a problem can be
formulated as follows: retrieve a log-concave function given its
area measure. Clearly, in view of (\ref{introduale1}) and
(\ref{introduale2}), the datum will consist of a single measure on
$\R^n$ or of a pair of measures (the first on $\R^n$ and the second
on $\sfe$), depending on whether we want to solve the problem in the
class ${\mathcal A}'$ or ${\mathcal A}''$, respectively. We
establish a uniqueness result for both these problems (see
Proposition \ref{uni}), and we find some necessary conditions for
the existence of a solution, which are quite similar to those afore
mentioned about the classic Minkowski problem (see Proposition
\ref{neccond}). However, differently from the case of convex bodies,
it turns out that such conditions are in general {\it not}
sufficient, as the analysis of the one dimensional case easily
shows. Thus, at this stage,  some substantial difference between the
geometric and the functional setting emerges, which deserves in our
opinion further investigation.

\medskip

\noindent{\bf Acknowledgements.} The authors wish to thank Shiri Artstein, Vitali Milman and Liram Rotem for
their suggestions of useful bibliographical references related to the subject of this paper.

\section{Preliminaries}

\subsection{Notation and background}\label{sec1}

We work in the $n$-dimensional Euclidean space $\R^n$, $n\ge1$,
endowed with the usual scalar product $\langle x ,  y \rangle $ and
norm $\|x\|$; we set $B _r := \{ x \in \R ^n \, :\, \|x\| \leq r \}$.

For $m \leq n$, we indicate by $\H^m$ the $m$--dimensional Hausdorff
measure; integration with respect to the Lebesgue measure $\H^n$ is
abbreviated by $dx$.

We denote by $\K ^n$ the class of convex bodies (compact convex
sets) in $\R^n$, and by $\K _0 ^n$ the subclass of convex bodies $K$
whose relative interior ${\rm int} (K)$ is nonempty. We indicate by
$V(K) = \H ^n (K)$ the $n$-dimensional volume of $K \in \K ^n$.

Given $K \in \K _0 ^n$, we denote by $\nu _K$ its Gauss map, by
$\sigma _K = (\nu_K) _\sharp (\H ^ {n-1} \res \partial K)$ its
surface area measure, and by $P (K) = \int_{S^{n-1} } \, d \sigma _K
= \H ^ {n-1} (\partial K)$ its perimeter. We say that $K$ is $\C
^2_+$ if its boundary $\partial K$ is of class $\C^2$ with strictly
positive Gaussian curvature.

For any $K \in \K ^n$, we adopt the standard notation $h _K$ for the
{\it support function} of $K$, defined by
$$
h _K (x):= \sup_{ y \in K} \langle x,  y \rangle \qquad \forall x \in \R ^n\ .
$$

We recall that the {\it polar body}  $K ^o$  of $K$ is given by
$$K ^o := \{ y \in \R ^n \ :\ \langle x, y \rangle \leq 1 \ \forall
x \in K \}\ ;$$  if $0 \in {\rm int} (K)$,  the support function of ${K}$ agrees with the {\it gauge function} of $K^o$, namely
$$h _{K } (x) = \rho _{K^o }(x) := \inf \{ t \geq 0 \ : \ x \in tK^o \}\ .$$

We denote by $I_K$ and $\chi _K$ the {\it indicatrix function} and
{\it characteristic function} of $K$, defined respectively by

$$I _K(x) := \begin{cases}
0 &  \hbox{ if } x \in K \\
+ \infty & \hbox{ if } x \not \in K\, ,
\end{cases} \qquad
\chi _K (x) := \begin{cases}
1 &  \hbox{ if } x \in K \\
0 & \hbox{ if } x \not \in K\ .
\end{cases}
$$

Let
$u : \R^n\to\R\cup\{+\infty\}$ be a convex function. We set
$$
\dom(u)=\{x\in\mathbb{R}^n\,:\, u(x)\in\mathbb{R}\}\,.
$$
By the convexity of $u$, $\dom(u)$ is a convex set. We say that $u$ is {\it proper} if $\dom (u) \neq \emptyset$.
We say that $u$ is {\it of class $\C^ 2_+$} if it is twice differentiable on ${\rm int} (\dom (u))$, with a positive definite Hessian matrix.
We denote by ${\rm epi} (u)$ the {\it epigraph} of $u$.

We recall that the {\it Fenchel conjugate} of $u$ is the convex function defined by:
$$
u^*(y)=\sup_{x\in\sR^n} \langle x , y \rangle -u(x) \qquad \forall y \in \R ^n.
$$

%We denote by $\ov u$ the {\it closure} of $u$, namely its
%lower semicontinuous envelope:
%$$\ov u (x) = \sup \big \{ g(x) \ \big  | \ g: \R^n \to \R \cup \{ + \infty \}\, , \ g \leq u\, , \ g \hbox { is lower semicontinuous} \big \}\ .$$
%We say that $u$ is closed if $\ov u = u$.

On the class of convex functions from $\R^n$ to
$\R\cup\{+\infty\}$, we consider the operation of {\it infimal convolution}, defined by
\begin{equation}\label{sommaL}
u \Box v (x) := \inf _{y\in \sR^n} \big \{ u (x-y) + v (y) \big \} \qquad \forall x \in \R ^n\ ,
\end{equation}
and the following {\it right scalar multiplication} by a nonnegative real number $\alpha$:
\begin{equation}\label{prodottoL} (u \alpha) (x) := \begin{cases}
\alpha u \left ( \frac{x}{\alpha} \right ) &  \hbox{ if } \alpha >0 \\
I _{\{0\}} & \hbox{ if } \alpha = 0
\end{cases} \qquad \forall x \in \R ^n\ .
\end{equation}

Notice that these  operations are convexity preserving, and  that
the function $I _{\{0 \}}$ acts as the identity element
in (\ref{sommaL}).

The proposition below gathers some elementary properties of the
Fenchel conjugate, in particular about its behaviour with respect to
the operations defined above. For the proof, we refer to \cite{Ro}.
\medskip
\begin{prop}\label{primolemma}
Let
$u : \R^n\to\R\cup\{+\infty\}$ be a convex function. Then:
\begin{itemize}
\item[(i)] it holds $u ^*(0) = - \inf (u)$; in particular, $\inf (u) > - \infty$ implies $u ^*$ proper;
\item[(ii)] if $u$ is proper,  then $u^*(y)>-\infty$ $\forall y\in\R^n$;
%\item[(iii)] if $u$ and $u ^*$ are proper, then $u ^ {**} = \ov u$;
\item[(iii)] $\dom (u \Box v) = \dom (u) + \dom (v)$;
\item[(iv)] $(u \Box v)^* = u ^* + v ^*$;
\item[(v)] $(u \alpha)^* = \alpha u ^*$.
\end{itemize}
\end{prop}

%In particular in (iii) the sum of sets is the usual vector or
%Minkowski addition: $A+B=\{a+b\,:\, a\in A\,,\,b\in B\}$.

Given a differentiable real valued function
$u$ on an open subset $C$ of $\R^n$, the {\it Legendre conjugate} of the pair
$(C, u)$ is defined to be the pair $(D, v)$, where $D$ is the image of $C$ through the gradient mapping $\nabla u$, and
$$v (y) = \langle \nabla u ^ {-1} (y), y \rangle - u \big ( \nabla u ^ {-1} (y) \big ) \qquad \forall y \in D\ .$$
Such definition is well posed whenever, for any $y \in D$,  the value of $\langle x, y \rangle - u (x)$  turns out to be independent from the choice of the point
$x \in \nabla u ^ {-1} (y)$.

Following \cite{Ro}, we say that a pair $(C, u)$ is  a {\it convex function of Legendre type} if:
\begin{itemize}
\item[(a)] $C$ is a nonempty open convex set;
\item[(b)] $u$ is differentiable and strictly convex on $C$;
\item[(c)] $\lim _i\|\nabla u (x_i )\| \to + \infty$ whenever $\{x_i\}\subset C$ is a sequence converging to some $x \in \partial C$.
\end{itemize}

Within the class of convex functions of Legendre type, Fenchel and Legendre conjugates may be identified according to Proposition below \cite[Theorem 26.5]{Ro}.

\begin{prop}\label{proplegendre}
Let $u:\R ^n \to \R \cup \{+\infty\}$ be a closed convex function, and set $C := {\rm int}(\dom (u))$, $C ^ * := {\rm int}(\dom (u ^*))$.
Then $(C, u)$ is a convex function of Legendre type if and only if $(C ^*, u ^*)$ is. In this case, $(C^*, u ^*)$ is the Legendre conjugate of $(C, u)$ (and conversely). Moreover, $\nabla u: C \to C ^*$ is a continous bijection, with $(\nabla u) ^ {-1} = \nabla u ^*$.
\end{prop}

\subsection{Functional setting}

Let us introduce the classes of functions we deal with throughout the paper.

\medskip

\begin{definition} {\rm We set:
$$\begin{array}{ll}
&\L:=\{u : \R^n\to\R\cup\{+\infty\}\ \big | \ \mbox{$u$ proper,
convex,
$\lim_{\|x\|\to+\infty}u(x)=+\infty$}\}\,,
\\ \noalign{\smallskip}
&\A: =\{f: \R^n\to\R\ \big |\ f=e^{-u}\,,\,u\in\L\}\,.
\end{array}
$$}
\end{definition}

Below, we give some examples and basic properties of functions in $\L$;  we show that, consequently,
the class of log-concave functions $\A$ can be endowed with an algebraic structure
which extends in a natural way the usual Minkowski structure on $\K ^n$.

\medskip
\begin{example}

{\rm (i) For any $K \in \K ^n$, the function $u = I _K$ belongs to
$\L$.  Notice that $u ^* =h _K$ belongs to $\L$ if and only if $0 \in {\rm int} (K)$,  which shows that the class $\L$ is not closed under Fenchel transform.

\smallskip
(ii) For any $K \in \K ^n$ with $0 \in {\rm int} (K)$, and any $p \in [1, + \infty)$, the function $u = \frac{1}{p} h _K ^ p$ belongs to $\L$.
In particular, for any $p \in [1, + \infty)$, the function $u (x) = \frac{1}{p} \|x \| ^ p$ belongs to $\L$. }
\end{example}

\medskip
\begin{lemma}\label{preli} Let $u \in \L$. Then
there exist constants $a$ and $b$, with $a>0$,
such that
\begin{equation}\label{claimab}
u(x)\ge a\|x\|+b\qquad\forall\,x\in\R^n\,.
\end{equation}
Moreover
$u ^*$ is proper, and satisfies
$u ^*(y) > - \infty$  $\forall y \in \R ^n$.
\end{lemma}
\proof
In order to show (\ref{claimab}),
assume first that $0\in\dom(u)$.  Let $r>0$ be such that $u(x)\ge1 +u(0)$ if $\|x\|\ge r$;
for $\|x\|\ge r$, the convexity of $u$ implies
$$
u(x)\ge u(0)+ \Big(u\big(\frac{rx}{\|x\|}\big)-u(0)\Big)\,\frac{\|x\|}{r}\ge u(0)+\frac{\|x\|}{r}
\ .
$$
Then, setting $m:=\inf (u)$, it holds
$$
u(x)\ge m-1+\frac{\|x\|}{r} \qquad \hbox{ for } \|x\| \geq r\ .
$$
Since the above inequality is verified for $\|x\|\le r$ as well, it holds in $\R^n$. This shows that
(\ref{claimab}) is satisfied by taking $a = r ^ {-1}$ and $b = (m-1)$. In the general case, since $u$ is proper, one can choose
$x_0\in\dom(u)$, and apply the above argument to the function $u(x-x_0)$, which yields
$$
u(x)\ge a\|x+x_0\|+b\ge a\|x\|+b-a\|x_0\|\ .
$$

The properties of $u ^*$ follow from Proposition \ref{primolemma} (i) and (ii).\qed

\medskip\bigskip

We now use Lemma \ref{preli} in order to prove that $\L$ is closed under the operations of infimal convolution and right scalar multiplication
defined in
(\ref{sommaL}) and (\ref{prodottoL}).

\medskip
\begin{prop}\label{Lchiuso}
Let $u,v\in\L$ and $\alpha,\beta\ge0$. Then
$(u \alpha) \Box (v \beta) \in \L$.
\end{prop}

\proof
From definition (\ref{prodottoL}) it is immediate that $(u \alpha) \in \L$ for any $u \in \L$ and $\alpha \geq 0$.
So we have just to show that $u \Box v$ belongs to $\L$ for any $u, v \in \L$.
Set for brevity $w: =u  \Box v $.
Clearly, $w$ is a convex function defined in $\R^n$. Let us prove that
$w$ takes values into $\R \cup \{ + \infty
\}$,  is proper and diverges as $\|x\| \to + \infty$.

By Proposition \ref{primolemma} (i) and (iv), we have
$$\inf w = - w ^* (0) = - u ^ * (0) - v ^ * (0) = \inf ( u) + \inf (v) \ .$$
Since $\inf (u), \inf (v) > - \infty$, we infer that $\inf (w) > - \infty$, which shows that $w$
takes values into $\R \cup \{ + \infty
\}$.

By Proposition \ref{primolemma} (iii), $\dom (w) = \dom (u) + \dom (v)$, hence the properness of both $u$ and $v$ implies the same property for $w$.

Let $u (x) \geq a \|x\| + b$ and $v (x) \geq a' \|x\| + b'$ according to Lemma \ref{preli}, and set $c:= \min \{ a, a'\}>0$, $d:= b + b'$. We have
$$\begin{array}
{ll}
w (x) & = \displaystyle{\inf _{y\in \sR^n} \big \{ u ( x-y) + v (y) \big \}} \\ \noalign{\smallskip}
&\displaystyle{\geq \inf _{y\in \sR^n} \big \{ a\| x-y\| + b  + a' \| y\| + b' \big \}} \\ \noalign{\smallskip}
&\displaystyle{\geq \inf _{y\in \sR^n} \big \{ c (\| x-y\|+  \| y\|) + d \big \} \geq c \|x\| + d\ .}
\end{array}
$$
In particular, this implies that $w$ diverges as $\|x\| \to + \infty$. \qed

\medskip\bigskip

We are now in a position to endow the class $\A$ with an addition and a multiplication by nonnegative scalars.
These operations are internal to $\A$ thanks to
Proposition \ref{Lchiuso}.

\medskip
\begin{definition}\label{defmink}
{\rm Let $f=e^{-u}$, $g=e^{-v}$ $\in \A$, and let $\alpha, \beta\ge0$. We define
\begin{equation}\label{defminkformula}
\alpha\cdot f\,\oplus\,\beta\cdot g=e^{-[(u \alpha) \Box (v
\beta)]}\ ,
\end{equation}
which in explicit form reads
$$(\alpha\cdot f\,\oplus\,\beta\cdot g )(x) := \sup _{y \in \sR ^n} \ f \left ( \frac{x-y}{\alpha} \right ) ^ {\alpha}
g \left ( \frac{y}{\beta} \right ) ^ {\beta}\ .$$}
\end{definition}

\medskip
\begin{remark}{\rm
In view of the identities
$$
\begin{array}{ll}
&
%\displaystyle
{u \Box v (x)  = \inf \big \{ \mu \ :\ (x, \mu) \in {\rm epi} (u) +
{\rm epi} (v) \big \}}\\  \noalign{\vskip -.2cm} &
%\displaystyle
{(u \alpha) (x)  = \inf \big \{ \mu \ :\ (x, \mu) \in \alpha \, {\rm
epi}(u)  \big \}}\ ,
\end{array}
$$
the functional operation in (\ref{defminkformula}) has the following
geometrical interpretation: it corresponds to the Minkowski
combination with coefficients $\alpha$ and $\beta$ of the epigraphs
of $u$ and $v$ (as subsets of $\R ^{n+1}$).
%In the general case, the geometrical characterization of the infimal
%convolution of $u $ and $v$ is the following: it is the largest
%extended real-valued function whose epigraph contains the Minkowski
%sum of the epigraphs of $u$ and $v$. $\bem$ $\enm$
}
\end{remark}

\medskip

Next Proposition shows that, when restricted to suitable subclasses
of $\A$, Definition \ref{defmink} allows to recover different
algebraic structures on convex bodies. Recall that (see \cite{Lu}), for a fixed $p
\in [1, + \infty)$, the $p$-sum of $K$ and $L$ with coefficients
$\alpha$ and $\beta$ is the convex body $\alpha \cdot _p K + _p
\beta \cdot _p L$ defined by the equality
$$
h ^p_{\alpha \cdot _p K + _p \beta \cdot _p L} = \alpha  h_K ^ p + \beta h
_L ^ p\ .
$$

\medskip
\begin{prop}\label{oldex}
Set
$$\begin{array}{ll}
&\L _1 := \big \{
 h _{K^o} \ :\ K \in \K ^n\ ,\ 0 \in {\rm int} (K)
\big \} \\
&\L _q := \big \{ \frac{1}{q} (h _{K^o}) ^q \ :\ K \in \K ^n\ ,\ 0
\in {\rm int} (K) \big \}\ , \quad q \in (1, + \infty)\ , \\
& \L _\infty := \big \{ I _K \ :\ K \in \K ^n \big \}\ .
\end{array}$$
The above subclasses of $\L$ are closed with respect to the
operations defined in $(\ref{sommaL})$ and $(\ref{prodottoL})$.

More precisely, for any $\alpha, \beta \geq 0$, and any $u, v$
belonging to the same class $\L _q$, it holds
$$
(u \alpha) \Box (v \beta)  =  \begin{cases}
 h _{ K ^o \cap L ^o } & \mbox{ if $q=1$, $u=h_{K^\circ}$, $v=h_{L^\circ}$,} \\ \noalign{\smallskip}
\frac{1}{p} (h _{(\alpha \cdot _p K + _p \beta \cdot _p L) ^ o} ) ^
p & \mbox{with $p:= \frac{q}{q-1}$, if $q \in (1, + \infty)$, $u=\frac{1}{q}(h_{K^\circ})^q$, $v=\frac{1}{q}(h_{L^\circ})^q$,}  \\
\noalign{\smallskip} I_{\alpha K+\beta L}& \mbox{if $q = \infty$, $u=I_K$, $v=I_L$.}\ .
\end{cases}$$

\end{prop}

\proof Let $u \in \L _q$. We
have
\begin{equation}\label{calcoloF}
u ^*=\begin{cases}
I _{K^o}  & \hbox { if }  q=1 \\ \noalign{\smallskip}
 \frac{1}{p} h _{K } ^ p & \hbox{ if } q \in (1, + \infty) \\ \noalign{\smallskip}
h _K & \hbox{ if } q= \infty\ .
\end{cases}
\end{equation}
In particular, in order to check the above expression of $u ^*$ in case $q \in (1, + \infty)$, one can apply with $\phi (s) = \frac{s^q}{q}$ the following identity
holding for every increasing convex function $\phi$
(see {\it e.g.} \cite{HiLe}):
$$(\phi (h_{K^o})) ^ * (x) = \inf _{t \geq 0} \{\phi ^*
(t) + t h _{K^o} ^
*\big (\frac{x}{t})\}\ ; $$
this yields
$$\Big (\frac{1}{q} (h _{K^o}) ^q \Big ) ^* (x) = \inf _{ \{t \geq 0 \ :\ x \in tK^o\}} \big \{ \frac{t ^p}{p} \big \}  = \frac{1}{p}  \rho _{K^o} ^ p (x)   = \frac{1}{p} {h _{K } ^p (x) } \ .$$

Now,  the statement of the Proposition  follows easily from the computation of $((u \alpha) \Box (v \beta)\big ) ^*$. Indeed, by Proposition \ref{primolemma} (iv)-(v), it holds  $((u \alpha) \Box (v \beta)\big ) ^* = \alpha u ^* + \beta v ^*$. According to
(\ref{calcoloF}), one has
$$\alpha u ^* + \beta v ^* = \begin{cases}
\alpha I _{K^o} +  \beta I _{L^o} = I _{K^o \cap L^o} = (h _{K^o \cap L^o}) ^* & \hbox{ if } q=1 \\
\noalign{\smallskip}
 \frac{1}{p}  \big [ \alpha h _K ^p+ \beta h _K
^p \big ] = \frac{1}{p}  \big [ h _{\alpha \cdot _p K + _p \beta
\cdot _p L }     ] ^p =  \big \{ \frac{1}{q}  \big [ h _{({\alpha
\cdot _p K + _p \beta \cdot _p L })^o }    ] ^q \big \} ^* & \hbox{ if } q\in ( 1, + \infty)
\\ \noalign{\smallskip }
\alpha h _K + \beta h _L = h _{\alpha K + \beta
L } = (I_{\alpha K+\beta L}) ^* & \hbox{ if } q= \infty \ .
\end{cases}
$$

\qed

\section{Differentiability of the total mass functional}

\begin{definition} {\rm We call {\it total mass functional} the
following integral
$$
J(f)=\int_{\sR^n} f(x)\,dx\qquad \forall f \in \A\ .
$$
}
\end{definition}

\medskip
\begin{remark}\label{varieJ}
{\rm
(i) The growth condition from below (\ref{claimab}) satisfied by functions in $\L$ ensures that
$
J(f)\in [0,+\infty)$ for every $f\in\A.
$

(ii) Clearly, when $f = \chi _K$, one has $J (f)= V(K)$.

(iii) If $f = e ^ {-u}$ is such that $J(f)=0$, then $f=0$
$\H^n$--a.e. in $\R^n$. This implies that the convex set $\dom(u)$
is Lebesgue negligible, and hence its dimension does not exceed
$(n-1)$. }
\end{remark}

\medskip
\begin{remark}\label{remPL}{\rm
By the Pr\'ekopa--Leindler inequality,
for every $f,g\in\A$ and for every $t\in[0,1]$, it holds
$$
J((1-t)\cdot f\,\oplus\,t\cdot g)\ge J(f)^{1-t}\,J(g)^t\ ,
$$
with equality sign if and only if there exists $x_0\in\R^n$ such
that $ g(x)=f\left({x-x_0}\right)\ \forall\, x\in\R^n\, $ (see
\cite{Dubuc,Gardner}). Consequently, for every fixed $f, g \in \A$,
the functions $t \mapsto \log J (f \oplus t \cdot g)$ and $t \mapsto
\log J \big ( (1-t) \cdot f \oplus t \cdot g \big )$  turn out to be
concave respectively on $[ 0 , + \infty)$ and on $[0,1]$. We shall
repeatedly exploit this concavity property in the sequel.}
\end{remark}

We are going to study the first variation of the total mass
functional, with respect to the algebraic structure introduced in
Definition \ref{defmink}.

\medskip
\begin{definition}{\rm Let $f,g\in\A$. Whenever the following limit
exists
$$
\lim_{t\to0^+}\frac{J(f\oplus t\cdot g)-J(f)}{t}\,,
$$
we denote it by $\D$,  and we call it {\it the first variation of
$J$ at $f$ along $g$}.}
\end{definition}

\medskip
\begin{remark}\label{ex1} {\rm Let $f=\chi_K$ and $g=\chi_L$, with  $K, L \in \K ^n$.
In this case $J(f\oplus t \cdot g)=V(K+tL)$ is a polynomial in $t$;
its derivative at $t=0^+$ is equal to $n$ times the {\it mixed
volume} $V_1(K,L)$, and admits the integral representation
\begin{equation}\label{surfmeasure}
\frac{d}{dt} V(K+tL) _{|_{t= 0^+ }} = n V_1(K, L) = \int _{S
^{n-1} } h _L \, d \sigma _K\ .
\end{equation}
Notice in particular that $\delta J (\chi _K , \chi _L)$ is nonnegative
and finite, which is not always true in general for $\delta J (f, g)$ ({\it cf.}
the examples given in Remark \ref{ex2} below).}
\end{remark}

\medskip
Subsection \ref{existence} below is devoted to prove that $\D$ exists under the fairly weak
hypothesis that $J (f)$ is strictly positive.  Then in subsection \ref{examples}
we show the explicit expression of $\D$ in some relevant cases.

\subsection{Existence of the first variation}\label{existence}

\medskip\bigskip
\begin{teo}\label{diffe} Let $f,g\in\A$, and assume that $J(f)>0$. Then
$J$ is differentiable at $f$ along $g$, and it holds
\begin{equation}\label{estid}
\D \in[- k, + \infty]\,,
\end{equation}
being $ k:= [ \inf( -\log g) ]_+\, J (f) $.  In dimension $n = 1$,
the same conclusions continue to hold also when $J (f) = 0$.
\end{teo}

\medskip
\begin{remark} {\rm We point out that the assumption $J (f)>0$ it
somehow technical; we believe that, when $J (f)=0$, Theorem
\ref{diffe} is likely true not only in dimension $n=1$ but also in
higher dimensions (as it is suggested by the fact that the mixed
volume $V_1(K, L)$ exists regardless the condition
$V(K)>0$). }\end{remark}

\begin{remark}\label{ex2} {\rm Estimate
(\ref{estid}) cannot be improved, as the following examples show.

(i) Let $f = e ^ {-u} \in \A$ with $J(f)>0$, and $g=e^{-v}$, where
$v(0)=1$ and $v\equiv+\infty$ in $\R^n\setminus\{0\}$. Then $u \Box
(vt) (x) = u (x) + t$, which implies
$$\D= J (f) \cdot \lim _{t \to 0 ^+} \frac{ e ^ {-t} -1} {t} =   -J(f)<0\ .$$

(ii)  Let $K, L \in \K ^n$ with the origin in their interior, so
that $u=h_K, v=h_L\in\L$, and take $f=e^{-u}, g=e^{-v}$. Then $u
\Box (vt)  = h _{K \cap L}$ ({\it cf.} Proposition \ref{oldex}), and therefore
$$
\D=\lim_{t\to0^+}
\left[\frac{1}{t}
\,\int_{\sR^n}\big(e^{-h_{K\cap L}}-e^{-h_{L}}\big)\,dx \right]
 = \begin{cases} 0 & \mbox{ if  } L \subseteq K \\
+ \infty & \mbox{ otherwise }\ . \end{cases}
$$

}
\end{remark}

\bigskip
Prior to the proof of Theorem \ref{diffe}, we state a
preliminary lemma, which will be heavily exploited also in the next
section.

\begin{lemma}\label{mono} Let $f= e ^ {-u}, g = e ^ {-v} \in \A$. For $t \geq 0$, set $u _t = u \Box (vt)$ and $f_t = e ^ {-u _t}$.
Assume that $v (0) = 0$. Then,
%$\dom (u) \subseteq \dom (u _t)$ and,
for every fixed $x \in \R^n$, $u _t (x)$ and $f _t
(x)$ are respectively pointwise decreasing and increasing with
respect to $t$; in particular it holds
$$
u_1 (x) \leq u _t (x) \leq u (x) \quad \hbox{ and } \quad f (x) \leq
f _t (x) \leq f _1(x) \qquad \forall x \in \R^n \, ,\
\forall t \in [0,1]\ .
$$
%\item[(iii)] for every set $E \subset \subset \dom (u)$, the Lipschitz constant of $u _t$ on $E$ is bounded above by some positive constant  independent of $t$.
\end{lemma}
\proof
%The inclusion $\dom (u) \subseteq \dom (u _t)$ is an
%immediate consequence  of the equality $\dom (u _t) = \dom (u) + t
%\dom (v)$ and the assumption $0 \in \dom (v)$.
Given $t \geq 0$ and $\delta >0$, let us show that $u _{t + \delta} \leq u _t$, {\it i.e.}
$$u \Box \big  (v(t+ \delta) \big ) \leq u \Box (vt)\ .$$

If $t=0$, the above inequality reduces to $u \Box (v \delta) \leq u$. This is readily checked:
recalling definitions (\ref{sommaL}) and (\ref{prodottoL}), from the assumption $v(0) = 0$ we deduce
$$u \Box (v \delta) (x) = \inf _{y \in \sR ^n} \big \{ u (x-y) + \delta v \big (\frac{y}{ \delta} \big ) \big \} \leq u (x)
\qquad \forall x \in \R^n \ .
$$

If $t>0$, for every $x \in \R^n$ we have
$$\begin{array}{ll} u \Box \big  (v(t+ \delta) \big )  (x) &\displaystyle{= \inf _ {\xi \in \sR^n}\Big \{  u (x-\xi) + ( t + \delta) v \big ( \frac{\xi}{t+\delta} \big ) \Big \} } \\
\noalign{\smallskip}
&\displaystyle{= \inf _ {\xi \in \sR^n} \Big \{ u (x-\xi) + \inf _{y \in \sR ^n}  \Big [t v \big ( \frac{\xi-y}{t} \big ) + \delta v \big ( \frac{y}{\delta} \big )\Big ]  \Big \} } \\
\noalign{\smallskip}
&\displaystyle{= \inf _ {y, z \in \sR^n} \Big \{ u (x-y-z) +   t v \big ( \frac{z}{t} \big ) + \delta v \big ( \frac{y}{\delta} \big )  \Big \} } \\
\noalign{\smallskip}
&\displaystyle{=  \big (u \Box \big  (vt) \big ) \Box (v \delta)   (x) \leq u \Box (v t) (x)   }. \\
\end{array}$$

Thus $u_t$ is monotone decreasing with respect to $t$, which immediately implies that $f_t = e ^ {- u _t}$ is monotone increasing.
\qed
\bigskip

{\it Proof of Theorem \ref{diffe}}.
We set
\begin{equation}\label{def1}
u := - \log f \, ,\qquad  v: = - \log g \, , \qquad  f_t:= f \oplus t
\cdot g\, , \end{equation} and
 \begin{equation}\label{def2}d:= v (0)\, , \qquad \tilde v (x) := v
(x) - d\, , \qquad  \tilde g(x) := e ^ {- \tilde v(x)}\, , \qquad \tilde f _t:= f \oplus t \cdot \tilde g \ .\end{equation}
Up to a translation of coordinates, we may also assume without loss of generality that $\inf (v) = v (0)$.

Since by construction
$\tilde v(0) = 0$, by Lemma \ref{mono} for every $x \in \R^n$ there exists $\tilde f(x):= \lim _{ t \to 0 ^+} \tilde f _t (x)$ and it holds $\tilde f(x) \geq f(x)$. Moreover, by monotone convergence, we have
$\lim _{t \to 0 ^ +} J(\tilde f_t) = J (\tilde f)$.

Since $f_t (x) = e ^ {-dt} \tilde f _t (x)$, we have
\begin{equation}\label{2.3}
\frac{J(f_t)-J(f)}{t} =J(f)\,\frac{e^{-dt}-1}{t} +
e^{-dt}\,\frac{J(\tilde f_t)-J(f)}{t}\, .
\end{equation}
Let us consider separately the two cases $J(\tilde f)>J(f)$ and $J(\tilde f)=J(f)$.

If $J(\tilde f)>J(f)$, then
$$
\lim_{t\to0^+}\frac{J(f_t)-J(f)}{t}= \lim_{t\to0^+}\frac{J(\tilde
f_t)-J(f)}{t}=+\infty\ ,
$$
and the thesis of the theorem holds true.

If $J (\tilde f) = J (f)$, we further distinguish the following two subcases:
$$
\exists t_0 >0 \ :\ J (\tilde f_{t_0}) = J (f) \qquad \hbox{ or } \qquad
J(\tilde f_t)>J(f) \ \forall t>0 \ .
$$In the former subcase, since by Lemma \ref{mono} $J (\tilde f_t)$ is a monotone increasing function of $t$, necessarily it holds $J (\tilde f_{t_0}) = J (f) $ for every $t \in [0, t_0]$.
Hence the
second addendum in the r.h.s.\ of (\ref{2.3})  is infinitesimal, so that
$$
\lim_{t \to 0 ^+} \frac{J(f_t)-J(f)}{t}
%=\frac{e^{-ct}J(\tilde f_t)-J(f)}{t}\\
=-d J(f)\,
$$
and the thesis
of the theorem holds true.

In the latter subcase, we can write
\begin{equation}\label{i}
\frac{J(\tilde f_t)-J(f)}{t}=
\frac{\log(J(\tilde f_t))-\log(J(f))}{t}\cdot
\frac{J(\tilde f_t)-J(f)}{\log(J(\tilde f_t))-\log(J(f))}\,.
\end{equation}
Since $\log(J(\tilde f_t))$
is an increasing concave function of $t$ (respectively by Lemma \ref{mono} and by the Pr\'ekopa--Leindler inequality, {\it cf.} Remark \ref{varieJ}), \begin{equation}\label{ii}
\exists \, \lim_{t\to0^+}\frac{\log(J(\tilde f_t))-\log(J(f))}{t}\in[0,+\infty]\,.
\end{equation}
On the other hand,
\begin{equation}\label{iii}
\exists\, \lim_{t\to0^+}\frac{J(\tilde f_t)-J(f)}{\log(J(\tilde f_t))-\log(J(f))}=J(f)>0\\ .
\end{equation}
From (\ref{i}), (\ref{ii}), and (\ref{iii}), we infer that
\begin{equation}\label{iv}
\exists \, \lim_{t\to0^+}\frac{J(\tilde
f_t)-J(f)}{t} \in [0, + \infty]\, .
\end{equation}
Combining (\ref{2.3}) and (\ref{iv}), we deduce that
\begin{equation}\label{v}
\exists \, \lim_{t\to0^+}\frac{J(f_t)-J(f)}{t}
\in [ - \max \{ d, 0 \} J (f), + \infty ]\,.
\end{equation}

Finally, let us show that in the one--dimensional case $\D$ exists
also when $J(f)=0$. We keep definitions (\ref{def1}) and (\ref{def2}). Since by assumption  $\dom (u)$ is a Lebesgue negligible convex set, it consists of
exactly one point $x_0$. Then
$$
u\square (\tilde vt)(x)=u(x_0)+t \tilde v\big (\frac {x -x_0}{ t}\big )\qquad \forall x \in \R\, , \ \forall t>0\,.
$$
Hence
\begin{equation}\label{vi}
\lim_{t\to0^+}\frac{J(\tilde f_t)-J(f)}{t}=
%\lim_{t\to0^+}\frac{e^{-u(x_0)}}{t}\int_{\sR} e^{-t\tilde v(x - x_0/t)}\, dx=
\lim_{t\to0^+}e^{-u(x_0)}\int_{\sR} e^{-t\tilde v(x- x_0)} \, dx=
e^{-u(x_0)}\H^1(\dom(v)) \in [0, + \infty]\,,
\end{equation}
where the last equality holds true by monotone convergence. Combining (\ref{2.3}) and (\ref{vi}), we see that (\ref{v}) remains true.
\qed

\medskip
\bigskip

\subsection{Computation of the first variation in some special cases}\label{examples}

Firstly, we analyze
the case $f=g$, and we show
that $\delta J (f, f)$ admits a very simple representation in terms
of the mass
 and the entropy of $f$, intended according to the definition below ({\it cf.} \cite{Ledoux}).

\begin{definition}\label{defent} {\rm For every $f \in \A$ with $J (f) >0$, we call {\it entropy} of $f$ the following quantity:}
$$
\Ent(f) = \int_{\sR^n} f \log f \, dx - J (f) \log J (f) \qquad
\forall f \in \A\ .
$$
\end{definition}
\medskip
\begin{prop}\label{ff}
For every $f \in \A$ with $J (f) >0$, it holds $\Ent(f) \in ( - \infty, + \infty)$ and
\begin{equation}\label{forment}
\delta J (f, f) =nJ(f)   +\int_{\sR^n} f \log f \, dx = \big ( n +
\log J (f) \big )J (f) + \Ent(f)\,.
\end{equation}
\end{prop}

\proof Since $J (f) \in (0, + \infty)$ for every $f \in \A$, to prove the finiteness of $\Ent (f)$ we have just to show that
$$\int_{\sR^n} f \log f \, dx \in ( - \infty , + \infty)\ .$$
We set $u : = - \log f$ and $\Omega:=\{x\in\R^n\,:\,u(x)\le0\} $
(which is possibly an empty set). It holds
$$
\int_{\Omega} f\log f \,dx= - \int_{\Omega} fu\,dx < - \inf _{\Omega} (u) \int _\Omega f < + \infty\ ,
$$
where in the last inequality we have used the boundedness of $u$ from below on $\Omega$ and the finiteness of $J (f)$.
On the other hand, we have
$$
\int_{\sR^n \setminus \Omega }  f\log f\,dx=- \int_{\sR^n\setminus\Omega} fu\,dx \geq
- m\, \int_{\sR^n\setminus\Omega} e^{-u(x)/2} \,dx > - \infty\ ,
$$
where we have used the elementary inequality ${t}{e^{-t/2}}\le m:=
2/e $ holding for every $t \in \R_+$ and Lemma \ref{preli}. So we have
$J(f \log f) \in ( - \infty, + \infty)$.

In order to prove the representation formula (\ref{forment}), assume first that $u\ge0$. Since
$
u\square (ut)= u(1+t)\ ,
$
we have
\begin{eqnarray*}
\frac{J(f\oplus t\cdot f)-J(f)}{t}&=& \frac
1t\left[(1+t)^n\int_{\sR^n}e^{-(1+t)u}\,dx-
\int_{\sR^n}e^{-u}\,dx\right]\\
\\
&=&\left[ \frac{(1+t)^n-1}{t} \right]\int_{\sR^n}e^{-(1+t)u}\,dx +
 \int_{\sR^n}e^{-u}\left( \frac{e^{-tu}-1}{t}
\right)\,dx \,.
\end{eqnarray*}
Now (\ref{forment}) follows by passing to the limit as $t \to 0 ^+$ (notice indeed that by the assumption $u \geq 0$ one can apply the monotone convergence theorem).

In the general case when the assumption
$u\ge0$ is removed, we consider the function $\tilde f=e^{-\tilde
u}$, where $\tilde u=u+c$ and $c=-\inf (u)$. One can easily check
that  $u\square (ut)=-c (1+t)+\tilde u\square (\tilde u t)$ and
consequently $J(f\oplus t \cdot f)=e^{c(1+t)}J(\tilde f\oplus t
\cdot \tilde f)$. As $\tilde u\ge 0$, we know that $\delta J (\tilde
f, \tilde f)$ exists and it is finite, so the same is true for
$\delta J (f, f)$. Moreover,
$$
\delta J (f, f) =ce^{c}J(\tilde f)+e^c \delta J (\tilde f, \tilde
f)=cJ(f) +e^c\left[ nJ(\tilde f)-\int_{\sR^n}e^{-(u+c)}(u+c)\,dx
\right]=n J(f)+ \int_{\sR^n} f \log f \, dx\,.
$$

\qed

\bigskip
\begin{remark}\label{sx}{\rm By inspection of the above proof, one can readily check that also the left derivative
$$\lim_{t \to 0 ^-} \frac{J (f \oplus t \cdot f) - J (f)} {t}$$
exists and agrees with $\delta J (f, f)$.}
\end{remark}

\bigskip
Next we consider the case when $f$ and $g$ belong to the class $\A_q$ introduced in Proposition \ref{oldex},  and we show that
$\delta J(f, g)$ can be written explicitly in integral form, by using the representation formula for $p$-mixed volumes
given in \cite{Lu}.

\begin{prop}{
Let $q \in ( 1, + \infty)$, and let $p:= q/ (q-1)$. Let $K, L \in \K ^n$ with the origin in their interior,
let $u:= \frac{1}{q} (h _{K^o}) ^q$, $v:= \frac{1}{q} (h _{L^o}) ^q$, and $f:= e ^ {-u}$, $g:= e ^ {-v}$.
There exists a positive constant $c = c (n, q)$ such that
\begin{eqnarray} & J (f) =
c(q, n)  V (K) & \label{massa}
\\ \noalign{\medskip}
& \displaystyle{\delta J (f, g) =
  \frac{c(q, n)}{n} \int _{S ^ {n-1} } h _L (\xi) ^p (h
_K(\xi))  ^ {1-p} \, d \sigma _K (\xi)}\ .  & \label{derivata}
\end{eqnarray}
}
\end{prop}

\proof   We set for brevity $a(t) = t ^ p/p$, so that  $a ^
* (t) = t ^ q/q$. We have:
$$\begin{array}{ll}J (f) = \displaystyle{\int _{\sR ^n} e ^ {-a^*(h _{K ^o})} \, dx} & = \displaystyle{\int _0 ^ 1 \H ^n
\big ( \big \{ x \ :\ e ^ {-a^*(h _{K ^o})(x)} > t \big \} \big ) \,
dt }
\\ \noalign{\medskip} & \displaystyle{= \int _0 ^ 1 \H ^n \big
( \big \{ x \ :\ h _{K ^o}(x) < (a^*) ^ {-1} (-\log t) \big \} \big
) \, dt}
\\ \noalign{\medskip} & \displaystyle{= \int _0 ^ 1 \H ^n \Big
( \Big \{ x \ :\ h _{K ^o}\Big ( \frac{x}{(a^*) ^ {-1} (-\log t)}
\Big ) < 1 \Big \} \Big ) \, dt}
\\ \noalign{\medskip} & \displaystyle{= \int _0 ^ 1 \big ((a^*) ^ {-1} (-\log t) \big ) ^n \, \H ^n \big
( \big \{ y \ :\ h _{K ^o}\big ( {y} \big ) < 1 \big \} \big ) \,
dt} \\ \noalign{\medskip} & \displaystyle{= \Big \{\int _0 ^ 1 \big ((a^*) ^ {-1} (-\log t) \big ) ^n \, dt  \Big \} \, V (K )\ ,}
\end{array}
$$
which proves (\ref{massa}) with $c(q, n):= \int _0 ^ 1 \big ((a^*) ^ {-1} (-\log t) \big ) ^n \, dt$.

Now we recall from Proposition \ref{oldex} that
$$f \oplus t \cdot g = e ^ {-\frac{1}{q} (h _{(K + _p t \cdot _p L ) ^o} ) ^ q}\ ,$$
which combined with (\ref{massa}) implies
$$\D = c (q, n) \lim _{t \to 0 ^ +} \frac{V ( K +_pt \cdot _p L ) - V
(K)} {t}\ . $$ Then (\ref{derivata}) follows from  the representation
formula for $p$-mixed volumes given in \cite[(IIIp)]{Lu}. \qed

\section{Integral representation of the first variation}\label{sec2}

In view of the examples in Section \ref{examples}, it is natural to
ask whether  $\delta J (f, g)$ admits more in general some kind of integral
representation. In this section we show that this is true when both $f$ and
$g$ belong to suitable subclasses of $\A$.

Let us begin by introducing the measures which
intervene in the representation formulae for $\D$. Such measures can
be viewed as the ``first variation'' of $J$ in the class of log-concave functions, since they play for $f$ the same role as the surface area measure for the volume in Convex Geometry. This fact emerges in a clear way
by comparing the first variation of volume in (\ref{surfmeasure}) with Theorems \ref{A'} and
\ref{A''} below.

\begin{definition}\label{defmusigma}
{\rm Let $f = e ^ {-u}\in \A$.  We set $\mu (f)$ the Borel measure
on $\R ^n$ defined by
$$\mu (f):= (\nabla u) _\sharp ( f \H ^n) \ .$$
When $\dom (u)= K \in \K ^n$, we also set $\sigma (f)$ the  Borel
measure on $S ^ {n-1}$ defined by
$$\sigma (f) := (\nu _K )_\sharp ( f \H ^ {n-1} \res \partial K)\ .$$}
\end{definition}

Next, we define the subclasses of $\A$ where our integral
representation formualae are settled.
\begin{definition}\label{subclasses} {\rm We set $\A', \A''$ the subclasses of $\A$ given by functions $f$ such that
$u  = - \log f $ belongs respectively  to
$$
\begin{array}{ll}
&
 \L ' := \Big \{ u \in \L \ :\ \dom (u)= \R^n
,\quad
 u \in \C ^ 2 _+ (\R ^n) \ ,\quad  \lim \limits_{\|x\| \to + \infty} \frac{u (x)} {\|x\|} = +
\infty \Big \}
\\  \noalign{\medskip}
& \L''  := \Big \{ u \in \L\ :\  \dom (u) = K \in \K ^n \cap \C ^2 _+, \quad u
\in \C^2_+ ({\rm int} (K) ) \cap \C ^0 (K)\ ,\quad \lim \limits_{x
\to
\partial K} \| \nabla u (x) \| = + \infty \Big \}
\ .
\end{array}
$$}
\end{definition}

\begin{remark}\label{legendre}
{\rm Notice that, for any $u$ belonging to $\L'$ or  $\L''$,  $({\rm int}(\dom (u)), u)$ is a convex function of Legendre type, and $u$ is {\it cofinite}, {\it i.e.}\ the domain of its Fenchel conjugate is the whole $\R^n$.} \end{remark}

Finally, we introduce the concept of admissible perturbation.
\begin{definition}\label{defH} {\rm We say that $g= e ^ {-v}$ is an {\it admissible perturbation} for
$f= e ^ {-u}$ if
\begin{equation}\label{H}
\exists \, c>0 \ :\ \varphi - c \psi \hbox{ is convex\ , \ where } \varphi= u ^* \hbox{ and } \psi = v ^*\ .
\end{equation}}
\end{definition}

Our integral representation results read as follows.

\begin{teo}\label{A'} Let $f,g \in \A'$, and assume that $g$ is an admissible perturbation for $f$.
Then $\D$ is finite and is given by
\begin{equation}\label{tesi'}
\D = \int _{\sR ^n} \psi \,d \mu (f)\ ,
\end{equation}
where $\psi = v ^*$.
\end{teo}

\begin{teo}\label{A''} Let $f,g \in \A''$, and assume that $g$ is an admissible perturbation for $f$.
Then $\D$ is finite and is given by
\begin{equation}\label{tesi''}
\D = \int _{\sR ^n} \psi \,d \mu (f)  + \int _{ S ^ {n-1}} h _L \, d
\sigma (f)\ ,
\end{equation}
where $\psi = v ^*$ and $L= \dom (v)$.
\end{teo}

\begin{remark}\label{remdim1} {\rm For $n=1$, $(\ref{tesi'})$ and $(\ref{tesi''})$ continue to
hold, possibly as an equality $+ \infty = + \infty$, if the
assumption that $g$ is an admissible perturbation for $f$ is
removed (see the Appendix for a proof).}
\end{remark}
\begin{remark}\label{comments}{\rm
Under the assumptions of Theorem \ref{A'} or Theorem \ref{A''},
by using the definition of push-forward measure and the change of variables $\nabla u
(x) = y$, one obtains
$$\int _{\sR ^n} \psi \,d \mu (f)= \int _{\dom (u)} \!\!\! \psi(\nabla u (x)) \, f(x) \, dx
=\int _{\sR ^n} \psi(y) \, e ^ { -\langle y,  \nabla \varphi (y) \rangle + \varphi (y)} \, \det (\nabla ^ 2 \varphi (y)) \, dy\ .
$$
Similarly, under the assumptions of Theorem \ref{A''}, it holds
$$
\int _{ S ^ {n-1}}\!\!\!\! h _L \, d \sigma (f) \!=\!
 \int _{
\partial K} \!\!\!\!h _L(\nu _K (x)) \, f (x)  \, d  \H ^ {n-1} (x) \!=\! \int _{ S ^ {n-1}}\!\!\!\! h _L(\xi)\,  f (\nu _K ^ {-1} (\xi)) \, \det (\nabla
\nu _K ^ {-1} (\xi)) \, d \H ^ {n-1} (\xi)\ .
$$}
\end{remark}

\bigskip

The proof of Theorems \ref{A'} and \ref{A''} is quite delicate and
requires several preliminary lemmas, whose proof is postponed to the Appendix.

The first one establishes the closure of the two
subclasses of $\L$ introduced in Definition \ref{subclasses} with
respect to the operations of infimal convolution and right scalar
multiplication.

\begin{lemma}\label{chiusura}
Let $u$ and $v$ belong both to the same class $\L'$ or $\L''$ and, for any $t >0$, set $u _t := u \Box (vt)$.
Then $u _t$ belongs to the same class as $u$ and $v$.
%\
%\begin{itemize}
%\item[(i)] Let $u,v\in\L'$ and $\alpha,\beta\ge0$. Then $(u \alpha)
%\Box (v \beta) \in \L'$.
%\item[(ii)] Let $u,v\in\L''$ and $\alpha,\beta\ge0$. Then $(u \alpha) \Box
%(v \beta) \in \L''$.
%\end{itemize}
\end{lemma}

We now turn attention to the behaviour of the functions $u _t = u \Box (vt)$ with respect to the parameter $t$, more precisely regarding their pointwise
convergence as $t \to 0 ^+$ (Lemma \ref{convergenza}), and their differentiability in $t$ (Lemma \ref{puntuale}).

\begin{lemma}\label{convergenza}
Let $u$ and $v$ belong both to the same class $\L'$ or $\L''$ and, for any $t >0$, set $u _t := u \Box (vt)$.
Assume that
$v(0) = 0$. Then
\begin{itemize}
\item[(i)] $\forall x \in \dom (u)\, , \qquad \, \lim _{t \to 0 ^+} u _t (x) =
u (x)$;
\item[(ii)] $\forall E \subset \subset  \dom (u)\, , \quad \lim _{t \to 0 ^+} \nabla u_t (x) =
\nabla u \hbox{ uniformly on } E$.
\end{itemize}
\end{lemma}

The following result is a key point in the proof of Theorems \ref{A'} and \ref{A''}; it contains an
explicit expression of the pointwise derivative of $u\Box(ut)$ with respect to $t$.

\begin{lemma}\label{puntuale}
Let $u$ and $v$ belong both to the same class $\L'$ or $\L''$ and, for any $t >0$, let $u _t := u \Box (vt)$.
Then
$$\forall x \in {\rm int} (\dom (u _t))\, , \ \forall t >0\,,
\qquad \frac{d}{dt} u _t (x) = -\psi (\nabla u _t (x) \big )\,,
\hbox{ where } \psi := v ^*\ . \qquad \qquad \qquad  $$
\end{lemma}

Next lemma provides a summability property of the Fenchel conjugate of $u = - \log f$ with respect to the measure
$\mu (f)$ introduced in Definition \ref{defmusigma}.

\begin{lemma}\label{finitezza}
Let $f = e ^ {-u} \in \A$, with $\varphi = u ^* \geq 0$. Then
$\varphi \in L ^ 1 (d \mu (f))$, namely
$$\int_{\sR^n} \varphi (\nabla u(x) ) \, f(x) \, dx < + \infty\ .$$
\end{lemma}

Finally, when $u$, $v \in \L ''$,  we need an estimate for $u _t = u \Box (vt)$ which will be exploited to deal with the boundary term in  Theorem \ref{A''}.

\begin{lemma}\label{ultimo}
Let $u, v \in \L ''$ and, for any $t >0$, let $u _t = u \Box (vt)$. Set $K:= \dom (u)$, $L:= \dom (v)$, $v_{\max}:= \max_L v$, and $v_{\min} := \min _L v$.
Then, for every $x \in K + tL$, there exists $y = y (x, t) \in K \cap (x - tL)$ such that
$$tv_{\min} + u (y) \leq u _t (x) \leq tv_{\max} + u (y)\ .$$
\end{lemma}

\bigskip

{\it Proof of Theorems \ref{A'} and \ref{A''}}.$\ $

We assume that either the hypotheses of Theorem \ref{A'} or the hypoteses of
Theorem \ref{A''}  are satisfied.

Throughout the proof we set
$$f = e ^ {-u}\ , \qquad g = e ^ {-v}\ ,
\qquad \varphi = u ^*, \qquad \psi = v ^ *\ ,\qquad E = \dom (u) \ , \qquad F = \dom (v)\ ,$$
and, for every $t \geq 0$,
$$f _t = f \oplus t \cdot g\ , \qquad u _t = u \Box (vt)\ ,\qquad \varphi _t =
\varphi + t\psi\ , \qquad E_t = E+ tF \ .$$
Let us point out that, under the assumptions of Theorem \ref{A'}, we have $E = F = \R ^n$,
whereas, under the assumptions of Theorem \ref{A''}, $E$ and $F$ are are convex bodies that will be named respectively $K$ and $L$.

Further, we need to `localize' our total mass functional: for every measurable set $A \subseteq \R ^n$ and any function $h \in \A$, we set
$$J _A (h) := \int _A h \, dx \ .$$

For convenience, we divide the proof into several steps.

\bigskip
{\it Step 1. Decomposition}.

\smallskip
With the notation introduced above, we can write
$$J ( f_t) - J (f) = J _ E ( f_t) - J_E (f) + J _{E _t \setminus E} (f _t) .$$
We are going to prove the integral representation formulae
(\ref{tesi'}) and (\ref{tesi''}) by showing that:

-- under the assumptions of one among Theorems \ref{A'} and \ref{A''}, it holds
\begin{equation}\label{interno}
\lim _{t \to 0 ^+} \frac{J _ E ( f_t) - J_E (f)}{t} = \int _{\sR ^n} \psi \,d \mu (f)\ ;
\end{equation}

-- under the assumptions of Theorem \ref{A''},  it holds
\begin{equation}\label{bordo}
\lim _{t \to 0 ^+}  \frac{J _{E _t \setminus E} (f _t)}{t} = \int _{ S ^ {n-1}} h _L \, d
\sigma (f)\ .
\end{equation}

\bigskip
 {\it Step 2. Reduction to the case  $0 \in {\rm int} (F)$, $v(0) = 0$, $v\geq 0$, $\varphi \geq 0$, $\psi \geq 0$.
}

\smallskip
Assume that equalities (\ref{interno}) and (\ref{bordo}) hold true (respectively under the assumptions of Theorems \ref{A'} or \ref{A''}, and of Theorem \ref{A''}), when  all the conditions $0 \in {\rm int} (F)$, $v(0) = 0$, $v\geq 0$, $\varphi \geq 0$, $\psi \geq 0$ are satisfied.

In the general case,
up to a translation of coordinates (which does not affect $J$), we may assume that
 $\inf v = v (0)$. Since by assumption $v$ belongs to $\L'$ or $\L ''$,
its minimum is necessarily attained in the interior of its domain, so we have $0 \in {\rm int} (F)$.
 If $c:= u
(0)$ and $d:= v(0)$, we set
 $$\tilde u (x) := u (x) - c ,\quad \tilde v (x) := v
(x) - d,\quad \tilde \varphi (y) := ( \tilde u) ^* (y),\quad
\tilde \psi (y) := ( \tilde v) ^* (y)\ $$
and
$$\tilde f  = e ^ {- \tilde u }, \quad \tilde g  = e ^
{- \tilde v }, \quad \tilde f_t := \tilde f \oplus t \cdot \tilde g
\ .$$ By construction it holds $\dom (\tilde v) = F$, $\tilde v(0) = 0$, $\tilde v\geq 0$,
$\tilde \varphi \geq 0$, $\tilde \psi \geq 0$.  Then, taking also into
account that $\dom (\tilde u) = E$,  $\tilde
\psi(y) = \psi(y) + d$, and $\tilde f = e ^ {c}  f$, it holds
\begin{equation}\label{tinterno}
\lim _{t \to 0 ^+} \frac{J _ E ( \tilde f_t) - J_E (\tilde f)}{t} = \int _{\sR ^n} \tilde \psi \,d \mu (\tilde f) =
e^ c \int _{\sR ^n}  \psi \,d \mu ( f) + d e ^c J _E(f)
\end{equation}
and
\begin{equation}\label{tbordo}
\lim _{t \to 0 ^+}  \frac{J _{E _t \setminus E} (\tilde f _t)}{t} = \int _{ S ^ {n-1}} h _L \, d
\sigma (\tilde f) = e ^c \int _{ S ^ {n-1}} h _L \, d
\sigma (f) \ .
\end{equation}
Now, since
$$f \oplus t \cdot g = e ^ {-(c+dt)}
(\tilde f \oplus t \cdot \tilde g) \ ,$$ we
may compute the left hand sides of (\ref{interno}) and (\ref{bordo}) as derivatives of a product.

Using (\ref{tinterno}), we get
$$\lim _{t \to 0 ^+} \frac{J _ E ( f_t) - J_E (f)}{t} =  - d e ^ {-c} J _E (\tilde f) + e ^ {-c}
\Big [ e^ c \int _{\sR ^n}  \psi \,d \mu ( f) + d e ^c J _E(f)
 \Big ] = \int _{\sR ^n}  \psi \,d \mu ( f)\ .$$

Similarly, using (\ref{tbordo}), we get
$$
\lim _{t \to 0 ^+}  \frac{J _{E _t \setminus E} (f _t)}{t} =
e ^ {-c} \cdot e ^c \int _{ S ^ {n-1}} h _L \, d
\sigma (f)   = \int _{ S ^ {n-1}} h _L \, d
\sigma (f)   \ .$$

\bigskip
{\it Step 3. For every $t > 0$, it holds
\begin{equation}\label{suicompatti}
J _E ( f_t ) -  J_E  (f) = \int _0 ^ t \Psi (s) \, ds\ ,
\end{equation}
where
\begin{equation}\label{defpsi} \Psi (s) : =
\displaystyle{\int_E \psi \, d \mu ( f_s)} \qquad \forall \, s\geq 0\ .
\end{equation}
}

\smallskip
Let $t>0$ be fixed, and take $C \subset \subset E$.  Thanks to the reduction $0 \in {\rm int}( F)$ made in Step 2, we have $C \subset \subset E_t$. Then
by Lemma \ref{puntuale} it holds
\begin{equation}\label{cpunt}
\lim _{h \to 0}  \frac{f_{t+h}(x) - f _t (x)}{h}  = \psi (\nabla u _t(x))  \, f_t(x) \qquad \forall x \in C \ .
\end{equation}
Moreover, thanks to the reduction $v(0) = 0$ made in Step 2, we can
apply Lemma \ref{mono} and Lemma \ref{convergenza} (ii) to infer
that, for every $s \in [0,1]$, the nonnegative functions $ \psi
(\nabla u _s(x))  \, f_s(x)$ are bounded above on $C$ by some
continuous function independent of $s$. Then, by the pointwise
convergence in (\ref{cpunt}), Lagrange theorem, and dominated
convergence we infer
$$\lim _{h \to 0} \frac{J _C(f_{t+h}) - J _C(f _t) }{h} \, dx = \lim _{h \to 0} \int _C \frac{f_{t+h} - f _t }{h} \, dx  = \int _C \psi (\nabla u _t)  \, f_t \, dx \ .$$
So we have
$$
J _C ( f_t ) -  J _ C (f) = \int _0 ^ t \Big \{ \int _C \psi \, d \mu (f_s) \Big \}\, ds\ ,
$$
which implies (\ref{suicompatti}) by letting $C \uparrow E$.

\bigskip
{\it Step 4. The function $\Psi$ defined in $(\ref{defpsi})$  takes
finite values at every $s\geq 0$.}

\smallskip
Let $s>0$. By the reduction $\varphi \geq 0$ made in Step 2,
we have
$$s \Psi (s)  \leq \int_{\sR ^n} (\varphi + s \psi) \, d \mu (
f_s) = \int _{\sR ^n}  u _s ^ * (\nabla u _s) f _s \, dx < + \infty\
,$$ where the last inequality follows from Lemma \ref{finitezza}
(which applies thanks to the conditions $\varphi, \psi \geq 0$).

Let now $s=0$. Since by assumption $g$ is an admissible perturbation for $f$, by (\ref{H}) it holds
$$(\varphi  - c  \psi) (y) \geq (\varphi  - c \psi) (0) +
\langle y, \nabla \varphi (0) - c \nabla \psi (0) \rangle\ ,$$ so
that
$$\psi (y) \leq c_1 + c _2 \varphi (y) + c _3 \|y \|\ ,$$
with
$$c_1 :=  \psi (0) - c^{-1} \varphi (0)\ , \qquad c_2:= c^ {-1}\ , \qquad
c_3 :=c^{-1} \|\nabla \varphi (0) - c \nabla \psi (0)\|\ .$$
Therefore
$$\begin{array}{ll} \displaystyle{\int_{\sR^n} \psi (\nabla u(x) ) \, f(x) \, dx }& \displaystyle{\leq c_1 \int_{\sR^n} f(x) \, dx
+  c _2 \int _{\sR^n} \varphi (\nabla u (x)) \, f (x) \, dx + c_3
\int _{\sR^n} \|\nabla u (x)\| \, f (x) \, dx}  \\
\noalign{\smallskip}& \displaystyle{=: c_1 I _1 + c _2 I _2 + c_3 I
_3}\ . \end{array} $$ Let us show separately that each of the
integrals $I _j$, $j = 1,2,3$, is finite. As already noticed in
Remark \ref{varieJ} (i), the integral $I _1$ is finite for every $f
\in \A$. The integral $I _2$ is finite by Lemma \ref{finitezza}.
Finally, in order to estimate the integral $I _3$, we use the coarea
formula: if $m:= \max_{\sR ^n} f$ it holds
\begin{equation}\label{i3prima}
I _ 3 = \int _{\sR ^n } \|\nabla f \| \, dx = \int _0 ^m {\mathcal
H} ^ {n-1} \big ( \partial \{f \geq s \} \big ) \, ds\ .
\end{equation}
According to Lemma \ref{preli}, there exist constant $a, b$, with
$a>0$ such that
$$f (x) \leq g(x):= e ^{ - a |x\| - b}\ ,$$
which implies $\{f \geq s \} \subseteq \{g \geq s \}$, and in turn,
\begin{equation}\label{i3seconda}
{\mathcal H} ^ {n-1} \big ( \partial \{f \geq s \} \big ) \leq
{\mathcal H} ^ {n-1} \big ( \partial \{g \geq s \} \big ) = c(n)
\Big ( \frac{-\log s -b}{a}  \Big ) ^{n-1} \ .
\end{equation}
The finiteness of $I _3$ follows from (\ref{i3prima}) and
(\ref{i3seconda}).

\bigskip
{\it Step 5.  The function $\Psi$ defined in $(\ref{defpsi})$ is
continuous at every $s>0$, and it is continuous from the right at
$s=0$.}

Through the change of variable $\nabla u _s (x)= y$, we obtain
$$\Psi (s) = \int_E \psi (\nabla u _s(x)) \, f_s (x) \, dx =
\int_{\sR^n} h (s,y) \, dy\ ,$$ with
$$h (s, y) := \psi (y) e ^ {\varphi _s (y) - \langle y, \nabla
\varphi _s (y) \rangle } \det (\nabla ^ 2 \varphi _s) (y) \chi
_{Q_s(y)}\ , \qquad   Q_s:= \nabla u _s (E) \ .$$ We now use
the expansion
$$\det (\nabla ^ 2 \varphi _s ) = \det (\nabla ^ 2 \varphi + s
\nabla  \psi ) = \sum _{j=0} ^n s ^j D_j (\varphi, \psi)\ ,$$ where
the mixed determinants $D_i (\varphi, \psi)$ are nonnegative functions of $y$
independent of $s$. We infer that
\begin{equation}\label{scompo}
\Psi (s) = \sum _{j=0} ^n s^j\Psi _j (s)\ ,
\end{equation}
where
$$\Psi _j (s): = \int_{\sR^n} h _j (s, y) \, dy \, \qquad h _j (s, y):=
\psi (y) e ^ {\varphi _s  - \langle y, \nabla \varphi _s
\rangle } D_j (\varphi, \psi) \chi _{Q_s}(y)\ .$$

Let us prove the continuity of $\Psi$ at a fixed $s_0 >0$. In view of (\ref{scompo}) it is enough to show that, for
any fixed index $i \in \{ 0, 1, \dots, n \}$, the function $\Psi _i$
is continuous at $s_0$.

We begin by noticing that \begin{equation}\label{convcar}
\lim _{s \to s_0} \chi _{Q_s}(y) = \lim _{s \to s_0} \chi _{Q _{s_0}}(y) \qquad \forall y \in \R^n\ .
 \end{equation}
Indeed, when $E = F= \R^n$, (\ref{convcar}) is trivially true since $Q_s = \R^n$ for every $s\geq 0$.
Assume $E = K$ and $F =L$, with $K, L \in \K ^n$.  The reduction $0 \in {\rm int} (F)$ made in Step 2 ensures that
$K \subset \subset K _{s_0}$, and hence by Lemma \ref{convergenza} (ii), we know that
$\nabla u _s$ converge uniformly to $\nabla u _{s_0}$ on $K$. Therefore, the compact sets  $Q_s$ converge to $Q_{s_0}$ in Hausdorff distance,
 which implies that the characteristic functions $\chi _{Q_s}$ converge to $\chi _{Q_{s_0}}$ in $L ^ 1 (\R^n)$, which in turn implies (\ref{convcar}).

Using (\ref{convcar}), we deduce that we have the pointwise convergence
$$
\lim _{s\to s_0} h _i (s, y) = h _i (s_0, y) \qquad \forall y \in \R^n\ .
$$

We claim that, as a consequence,  $\Psi _i (s)$
tends to $\Psi _i(s_0)$ as $s \to s _0$ by dominated convergence.
Indeed, let us show that $h _i (s, y)$ are bounded from above by a function in $L ^ 1 (\R ^n)$ independent of $s$. By the reduction $v \geq 0$ made in Step 2, for any fixed
$y \in \R ^n$ the map
\begin{equation}\label{crumon}
s \mapsto e ^{\varphi _s (y) - \langle y, \nabla \varphi _s (y)
\rangle } \end{equation} is pointwise decreasing. Therefore, if we
fix $\ov s \in (0, s_0)$, for any $s \geq \ov s$ it holds
$$\begin{array}{ll}
h _i (s, y) & \displaystyle{\leq \psi  e ^ {\varphi _{ \ov s}  -
\langle y, \nabla
\varphi _{\ov s}  \rangle } D_i (\varphi, \psi) } \\
&  \displaystyle{= \frac{1}{\ov s^i}  \psi  e ^ {\varphi _{\ov s}  -
\langle y,
\nabla \varphi _{\ov s}  \rangle } \ov s ^ i D_i (\varphi, \psi)} \\
&  \displaystyle{\leq \frac{1}{\ov s^i}  \psi  e ^ {\varphi _{\ov s}
- \langle y, \nabla \varphi _{\ov s}  \rangle }  \sum _{j=0} ^n \ov
s ^ j D_j (\varphi, \psi)}
\\
&  \displaystyle{= \frac{1}{\ov s^i}  \psi  e ^ {\varphi _{\ov s}  -
\langle y, \nabla \varphi _{\ov s}  \rangle }  \det (\nabla ^2
\varphi _{\ov s})}
\\
&  \displaystyle{\leq \frac{1}{\ov s^{i+1}}  \varphi _{\ov s}  e ^
{\varphi _{\ov s}  - \langle y, \nabla \varphi _{\ov s}  \rangle }
\det (\nabla ^2 \varphi _{\ov s})\ , }\end{array}
$$
and the function in the last line belongs to $L ^ 1 (\R ^n)$ by
Lemma \ref{finitezza}.

Let us now prove the continuity from the right of $\Psi$ at $s=0$. To that
aim, in view of (\ref{scompo}) is is enough to show that
\begin{eqnarray}
& \lim _{s \to 0 ^ +} \Psi _0 (s) = \Psi (0) \hskip 3 cm & \label{I}
\\ & \limsup _{s \to 0 ^ +} \Psi _i (s) < + \infty \qquad \hbox{  $\forall \, i \in \{1, \dots, n \}$.} & \label{II}
\end{eqnarray}
To prove equality (\ref{I}),  we begin by noticing that,
as $s \to 0 ^+$, the sets $Q _s$ invade $\R^n$,
meaning
\begin{equation}\label{inclusion}
\forall r>0, \ \exists \ov s>0 \ :\  Q_s  \supseteq B_r  \qquad
\forall s \leq \ov s\ .
\end{equation}
Indeed, when $E = F= \R^n$, (\ref{inclusion}) is trivially true since $Q_s = \R^n$ for every $s\geq 0$.
Assume $E = K$ and $F =L$, with $K, L \in \K ^n$, and
let $r>0$ be fixed. We have
\begin{equation}\label{defC}
Q _ s = \nabla u _s (K) \supseteq \nabla u _s (C) \ , \quad \hbox{ with }
C:= \nabla u ^ {-1} (B _{2r})\ .
\end{equation}
Since $C \subset \subset K$ and $K \subset \subset K _s$ (the latter  thanks to the reduction $0 \in {\rm int} (L)$ made in Step 2),
by Lemma \ref{convergenza} (ii) we know that
$\nabla u _{s}$ converge uniformly to $\nabla u$ on $C$.
Therefore, the compact sets  $\nabla u _s (C)$ converge to $B _{2r}$ in Hausdorff distance, so that they contain $B _r$ for $s$ sufficiently small.
Combined with (\ref{defC}), this implies (\ref{inclusion}).

%Therefore, there exists $\ov s$ such that
%$$\|\nabla u _{s} (x) \| > r \qquad \forall x \in \partial C\ ,\ \forall s \leq \ov s\ ,$$
%which implies
%$$\partial C_s \cap B _r = \emptyset \qquad \forall  s\leq \ov s\ .
%$$
%Since we may also assume that $ 0 \in F_s$, the above condition
%implies that $F_s \supseteq B _r$ $\forall s \leq \ov s$. Then by
%construction it holds $Q _s \supseteq F _s \supseteq B _r$ $\forall
%s \leq \ov s$, and (\ref{inclusion}) is proved.
Using (\ref{inclusion}), we deduce that
we have the pointwise convergence
\begin{equation}\label{pointwise2}
\lim _{s\to 0} h _0 (s, y) = h _0 (0, y)\qquad \forall y \in \R
^n .
\end{equation}
Now, by the monotonicity of the map (\ref{crumon}), for
any $s \geq 0$ it holds
\begin{equation}\label{dominante}
h _0 (s, y) \leq h _0 (0, y)= \psi  e ^ {\varphi - \langle y, \nabla \varphi \rangle} \det (\nabla ^ 2 \varphi )
\ ,
\end{equation}
and the last expression is in $L ^ 1 (\R ^n)$ because we have
proved in Step 4 that $\Psi (0)$ is finite.

In view of (\ref{pointwise2}) and (\ref{dominante}),  (\ref{I}) holds true by dominated convergence.

To prove (\ref{II}) we notice that assumption (\ref{H}) implies
$\nabla ^ 2 \psi \leq c ^ {-1} \nabla ^ 2 \varphi$ and hence
$$D_i (\varphi, \psi) \leq D_i ( \varphi , c ^ {-1} \varphi)\ .$$
This, combined with the monotonicity of the map (\ref{crumon}),
implies
$$h _i (s, y) \leq \psi e ^ {\varphi - \langle y, \nabla \varphi
\rangle} D_i ( \varphi , c ^ {-1} \varphi)  = \psi e ^
{\varphi - \langle y, \nabla \varphi \rangle}  \gamma _i (c) \det
(\nabla ^ 2 \varphi )  \ ,
$$ where the coefficients $\gamma _i (c)$ depend only on $c$.  The last expression is in $L ^ 1(\R^n)$ again
by the finiteness of $\Psi (0)$, and (\ref{II}) follows.

\bigskip
{\it Step 6. Equality $(\ref{interno})$ holds.}

\smallskip
The equality (\ref{suicompatti}) proved in Step 3, together with the
finiteness and  continuity of $\Psi (s)$ for $s>0$ proved
respectively in Steps 4 and 5, gives
\begin{equation}\label{positives}
\Psi (s) = \frac{d}{dt} J
_E(f _t) |_{t=s} \qquad \forall \, s>0\ .
\end{equation}  Moreover, the continuity
from the right of $\Psi$ at $s=0$ proved in Step 5 implies
\begin{equation}\label{caso1} \lim
_{s \to 0 ^+} \Psi (s) = \Psi(0) =\int _{\sR ^n } \psi \, d \mu (f)\
.
\end{equation}
Therefore,
\begin{equation}\label{limder}
\lim _{t \to 0 ^+ } \frac{J _ E ( f_t) - J_E (f)}{t} =
\frac{d}{dt} J _E(f _t) |_{t=0 ^+}= \lim _{s \to 0 ^+} \frac{d}{dt}
J _E(f _t) |_{t=s} =  \lim _{s \to 0 ^+} \Psi (s)
%= \Psi(0)
= \int
_{\sR ^n } \psi \, d \mu (f)\ .
\end{equation}

\bigskip
{\it Step 7. Under the assumptions of Theorem \ref{A''}, equality $(\ref{bordo})$ holds.}

\smallskip
We define the map $m:S ^ { n-1} \times [0,t] \to K _t \setminus K$ by
$$m  (\xi, s) := \nu _{K_s} ^ {-1} (\xi) = \nu _K ^ {-1} (\xi) + s \nu _L ^ {-1} (\xi) \ .$$
By the area formula \cite[Section 3.1.5]{Giaquinta}, we have
\begin{equation}\label{area}
\int_{K _t \setminus K } f _t  =  \int _0 ^ t \int _{S ^{n-1} } f_t (m (\xi, s)) |\det J m (\xi, s)| \, d \H ^ {n-1} (\xi) \, ds\ .
\end{equation}

Let $(\xi, s) \in S ^ {n-1} \times [0,t]$ be fixed and let us compute  $|\det J m (\xi, s)|$.
We choose an orthonormal basis $\{e_1, \dots, e _n \}$ of the tangent space $\xi ^ \perp \times \R$ to $S ^ { n-1} \times [0,t]$
given by
$$e _i = (v_i, 0) \ \quad i = 1, \dots, n-1\ , \qquad e _n = (0, \dots, 0 , 1)\ ,$$
where $v_i$ are eigenvectors of the reverse Weingarten operator $\nabla \nu _{K_s} ^ {-1}(\xi)$.
Then, denoting by $\rho_i (\xi, s)$ the corresponding eigenvalues (namely the principal radii of curvature of $\partial K _s$ at
$\xi$), it holds
$$\partial _{e_i} m (\xi, s) = \rho _i (\xi, s) e_i  \ \quad i = 1, \dots, n-1\ , \qquad  \partial _{e_n} m (\xi, s) = \nu _L ^ {-1} (\xi)\ .$$
Hence
$$
|\det Jm (\xi, s)|  = \| \partial _{e_1} m (\xi, s) \wedge \dots \wedge \partial _{e_n} m (\xi, s) \|=
|\langle \xi, \nu _L ^ {-1} (\xi) \rangle| \cdot \prod _{i=1} ^ {n-1} \rho _i (\xi, s)
=  h _L (\xi)  \det \big ( \nabla \nu _{K_s} ^ {-1}  (\xi) \big )
\ ,
$$
where the last equality holds because, by the reduction $0 \in {\rm int} (L)$ made in Step 2, we have  $h _L \geq 0$.

Now we recall that the reverse Weingarten operator of $K _s$ is given by
$$\nabla \nu _{K_s} ^ {-1}=(h _{K _s}) _{ij} +  h _{K _s} \delta _{ij}\ ,$$
where indices $i$ and $j$ denote second order covariant derivation with respect to an orthonormal frame on $S ^ {n-1}$.
Therefore, as $h _{K _s} = h _K + s h _L$, we have
$$\nabla \nu _{K_s} ^ {-1}=(h _{K }) _{ij} +  h _{K} \delta _{ij} + s \big [ (h _{L }) _{ij} +  h _{L} \delta _{ij} \big ]\ ,$$
and hence
\begin{equation}\label{jacobian}
|\det  Jm (\xi, s)|  = h _ L (\xi) \Big [ \det \big ( \nabla \nu _{K} ^ {-1}  (\xi) \big )  +
\sum _{i = 1} ^ {n-1} \gamma_i (\xi) s ^ i \Big ] \ ,
\end{equation}
where $\gamma _i
(\xi)$ are continuous functions depending only on the curvatures of $\partial K$ and $\partial L$ at $\xi$.

Inserting (\ref{jacobian}) into (\ref{area}) and dividing by $t$ we obtain
\begin{equation}\label{penultima}
\begin{array}{ll} \displaystyle{\frac{1}{t} \int _{K _ t \setminus K} f _t \, dx
} & = \displaystyle{ \frac{1}{t} \int _0 ^ t \int _{ S ^ {n-1} }f_t (m (\xi, s))
h _L (\xi) \det \big ( \nabla \nu _{K} ^ {-1}  (\xi) \big )   \, d \H ^ {n-1} (\xi) \, ds   }
\\ \noalign{\medskip}
& + \displaystyle{ \sum _{i = 1} ^ {n-1}  \frac{1}{t} \int _0 ^ t  s ^ i\Big \{  \int _{ S ^ {n-1} } f_t (m (\xi, s)) h _L (\xi)
 \gamma_i (\xi)   \, d \H ^ {n-1} (\xi)\Big \} \, ds }\ . \end{array}
\end{equation}
We observe that
\begin{equation}\label{somma0} \lim _{t \to 0 ^+} \sum _{i = 1} ^ {n-1} \frac{1}{t} \int _0 ^ t  s ^ i\Big \{  \int _{ S ^ {n-1} } f_t (m (\xi, s)) h _L (\xi)
 \gamma_i (\xi)   \, d \H ^ {n-1} (\xi)\Big \}  = 0\ .
 \end{equation}
Indeed, for every $i = 1, \dots, n -1$, we have
$$\int _0 ^ t  s ^ i\Big \{  \int _{ S ^ {n-1} } f_t (m (\xi, s)) h _L (\xi)
 \gamma_i (\xi)   \, d \H ^ {n-1} (\xi)\Big \} \leq  (\sup _{\sR^n} f_1) \int _{ S ^ {n-1} } \!\! h _L
 \gamma_i  \, d \H ^ {n-1}
 \int _0 ^ t  s ^ i  \, ds\ ,
 $$
where we used the inequality
$f_t (x) \leq f_1(x)$ holding for every $x \in \R^n$ and every $t \in [0,1]$ by Lemma \ref{mono} (which applies thanks to the reduction $v(0)= 0$ made in Step 2).

By (\ref{penultima}) and (\ref{somma0}), to conclude the proof of Step 7 it is enough to show that
$$\lim _{t \to 0 ^+}
\frac{1}{t} \int _0 ^ t \int _{ S ^ {n-1} }f_t (m (\xi, s))
h _L (\xi) \det \big ( \nabla \nu _{K} ^ {-1}  (\xi) \big )   \, d \H ^ {n-1} (\xi) \, ds    = \int _{S ^ {n-1} } h _L \, d \sigma (f)\ ,$$
or equivalently
$$\lim _{t \to 0 ^+} \frac{1}{t} \int _0 ^ t \int _{ S ^ {n-1}} \big [ f _t (m (\xi, s) - f (m (\xi, 0)) \big ] h _L ( \xi)
\det \big ( \nabla \nu _{K} ^ {-1}  (\xi) \big )    \, d \H ^ {n-1} (\xi)= 0\ .$$
Such equality is clearly satisfied if
\begin{equation}\label{uniforme}
\lim _{t \to 0 ^+}  \sup _{s \in [0,t], \ \xi \in  S ^ {n-1} } \big | f _t (m (\xi, s)) - f (m (\xi, 0)) \big |  = 0\ .
\end{equation}

Let $s \in [0,t]$ and $\xi \in  S ^ {n-1}$. By Lemma \ref{ultimo} applied at the point $x := m (\xi, s ) \in \partial K _s \subset K _t$, there exists
$y \in K    \cap (x- tL)$ such that
$$tv_{\min} + u (y) \leq u _t (m (\xi, s)) \leq tv_{\max} + u (y)\ .$$
Hence
\begin{equation}\label{d1}
tv_{\min} + u (y) - u ( m (\xi, 0)) \leq u _t (m (\xi, s)) - u ( m (\xi, 0)) \leq tv_{\max} + u (y)- u ( m (\xi, 0)) \ .
\end{equation}
As $x \in m (\xi, 0)+  s L \subseteq  m (\xi, 0) + t L $, we have $m (\xi, 0) \in K    \cap (x- tL)$, and therefore
\begin{equation}\label{d2}
\| m (\xi, 0) - y \| \leq {\rm diam} \big ( K    \cap (x- tL) \big ) \leq t {\rm diam} (L)\ .
\end{equation}
By (\ref{d1}), (\ref{d2}) and the uniform continuity of $u$ on $K$, we infer that
$$
\lim _{t \to 0 ^+}  \sup _{s \in [0,t], \ \xi \in  S ^ {n-1} } \big | u _t (m (\xi, s)) - u (m (\xi, 0)) \big |  = 0\ ,
$$
and (\ref{uniforme}) follows.

\bigskip
{\it Step 8: Conclusion.}

\smallskip
Equalities (\ref{tesi'}) and (\ref{tesi''}) follow from Steps 1, 6, and 7.
Moreover,  the finiteness
of $\Psi (0)$ proved in Step 4 implies that $\int_{\sR^n} \psi \, d \mu (f) < + \infty$; on the other hand,
for any $K, L \in
\K ^n$, one has $\int_{S^{n-1}} h _L \, d \sigma (f) < + \infty$.
Therefore
$\delta J (f, g)$ is finite.

\medskip

\section{The functional form of Minkowski first
inequality}\label{secmink}

Minkowski first inequality states that
\begin{equation}\label{M1}
\lim_{t\to 0 ^+} \frac{V (K+tL)- V (K)} {t} =  n V_1(K, L)
\geq n V(K) ^{1-\frac{1}{n}} V(L) ^ {\frac{1}{n}} \qquad \forall K,
L \in \K _0^n\ ,\end{equation} with equality sign if and only if $K$
and $L$ are homothetic (see \cite[Theorem 6.2.1]{Schneider}).

The main result of this section provides a functional version of
such inequality:

\bigskip

\begin{teo}\label{teomink} Let $f, g \in \A$, and assume that $J (f)>0$. Then
\begin{equation}\label{mink1}
\delta J (f, g) \geq J (f) \big [\log J (g)
 +n \Big ] +  \Ent (f) \ ,
\end{equation}
with equality sign if and only if there exists $x_0\in \R ^n$ such
that $g(x) = f (x-x_0)$ $\forall x \in \R^n$.
\end{teo}

\begin{remark}{\rm We point out that, by choosing $f= \gamma _n$, Theorem \ref{teomink} allows to recover the
Urysohn-type inequality for the mean width of a log-concave function
proved in \cite[Proposition 3.2]{Klartag-Milman05} and \cite[Theorem
1.4]{Rotem2}. }
\end{remark}

\bigskip Before giving the proof of Theorem \ref{teomink}, let us
present a straightforward consequence of it, which  will be
exploited in Section \ref{towards} in order to get uniqueness in the
functional form of the Minkowski problem.
\begin{cor}\label{cormink}
Let $f_1, f_2 \in \A$, with $J (f_1) = J (f_2) >0$, and assume that
\begin{equation}\label{incroci}
\delta J (f_2, f_1) = \delta J (f_1, f_1)  \qquad \hbox{ and }
\qquad \delta J (f_1, f_2) = \delta J (f_2, f_2) \ .\end{equation}
Then there exists $x_0\in \R ^n$ such that $f_2(x) = f_1 (x-x_0)$
$\forall x \in \R^n$.
\end{cor}

\proof By the assumption  $J (f_i) >0$, we may apply inequality
(\ref{mink1}) (once with $f= f_1$ and $g= f_2$ and once with $f=
f_2$ and $g= f_1$); since $J (f_1) = J (f_2)$, we get
\begin{equation}\label{basica}
\delta J (f_1, f_2) \geq n J (f_1)  + \int_{\sR^n} f_1 \log f_1 \,
dx\qquad \hbox{ and } \qquad
 \delta J (f_2, f_1) \geq n J (f_2)  + \int_{\sR^n} f_2 \log f_2 \, dx  \ .
\end{equation}
By assumption (\ref{incroci}) and Proposition \ref{ff}, the two
inequalities in (\ref{basica}) may be rewritten respectively as
$$\delta J (f_2, f_2) \geq \delta J (f_1, f_1) \qquad \hbox{ and } \qquad
 \delta J (f_1, f_1) \geq \delta J (f_2, f_2) \ ,
$$
which implies that both hold with equality sign. Then $f_1$ and
$f_2$ are translates of each other by Theorem \ref{teomink}. \qed

\bigskip
We now turn to the proof of Theorem \ref{teomink}. We need the
following

\bigskip
\begin{lemma}\label{jerison}
Let $f, g \in \A$, and assume that $J (f)>0$. Then
$$\lim _{t \to 0 ^ +} \frac{J \big ((1-t) \cdot f \oplus t \cdot g \big ) - J (f)   }{t} = \delta J (f, g) - \delta J (f,
f)\ .$$
\end{lemma}
\proof For $t \in (0,1)$, we set
$$\alpha (t) :=  \frac{t}{1-t} \qquad \hbox{and} \qquad
f_{\alpha (t)} := f \oplus \alpha (t) \cdot g\ .$$ Let us write
\begin{equation}\label{zero}
\frac{J\big (( 1-t) \cdot f \oplus t \cdot g\big ) - J (f) }{t} =
 \frac{J \big
((1-t) \cdot f_{\alpha (t)}\big ) - J \big ( f_{\alpha (t)} \big
)}{t}+ \frac{J \big ( f_{\alpha (t)} \big ) - J (f)}{t} \ ,
\end{equation}
 and let us focus attention on the the first adddendum in the r.h.s.\ of (\ref{zero}).

For every fixed $t \in (0,1)$, we have
$$
\frac{J \big ((1-t) \cdot f_{\alpha (t)}\big ) - J \big ( f_{\alpha
(t)} \big )}{t} = \frac{\gamma_t (t) - \gamma _t (0)} {t}\ ,
$$
where the function $\gamma _t$ is defined by
$$\gamma _t(s) := J \big ( (1-s) \cdot f_{\alpha (t)} \big ) \qquad \forall
s \in (0,1)\ .$$ In view of Proposition \ref{ff} and Remark
\ref{sx}, the function $\gamma_t$ is differentiable on $(0,t)$, with
$$\begin{array}{ll}
\gamma _t' ( s) &= \displaystyle{- \delta J \big ( (1-  s)\cdot
f_{\alpha (t)}, (1- s)\cdot f_{\alpha (t)} \big ) } \\
\noalign{\smallskip} & = - n J \big ( (1-  s)\cdot f_{\alpha (t)}
\big ) - \displaystyle{\int_{\sR^n} (1- s) \cdot f_{\alpha (t)} \log
\big ( (1- s) \cdot f_{\alpha (t)} \big ) \, dx}
\\ \noalign{\smallskip} & = (1 -  s) ^ n \big [ - n J
\big ( f_{\alpha (t)} ^ {1-  s} \big ) - \displaystyle{\int_{\sR^n}
f_{\alpha (t)} ^ {1-  s}  \log \big (f_{\alpha (t)} ^ {1-  s} \big )
\, dx \big ]} \ .
\end{array}
$$ Then, for every fixed $t \in (0,1)$, we can apply Lagrange theorem to infer that there exists $\ov s \in (0, t)$ such
that
\begin{equation}\label{tre}
\frac{J \big ((1-t) \cdot f_{\alpha (t)}\big ) - J \big ( f_{\alpha
(t)} \big )}{t}  = \gamma _t' (\ov s)  = (1 - \ov s) ^ n \big [ - n
J \big ( f_{\alpha (t)} ^ {1- \ov s} \big ) - \int_{\sR^n} f_{\alpha
(t)} ^ {1- \ov s}  \log \big (f_{\alpha (t)} ^ {1- \ov s}  \big ) \,
dx \big ] \ .
\end{equation}

We are now ready to pass to the limit as $t \to 0 ^+$ in the r.h.s.\
of (\ref{zero}).

Concerning the first addendum, assume for a moment that the function
$v:= -\log g$ satisfies the condition $v(0) = 0$. In this case, by
Lemma \ref{mono}, as $t \to 0 ^+$ the functions $f _{\alpha (t)}(x)$
converge increasingly to some  pointwise limit $\tilde f (x)$ (which
is bounded above and below by some functions in $\A$). Then, by
monotone convergence, taking also into account that $\ov s \to 0^+$
as $t \to 0 ^ +$, we infer from (\ref{tre}) that
\begin{equation}\label{due}
\lim _{t \to 0 ^ +} \frac{J \big ((1-t) \cdot f_{\alpha (t)}\big ) -
J \big ( f_{\alpha (t)} \big )}{t} =
   - n
J \big ( \tilde f  \big ) - \int_{\sR^n} \tilde f \log \tilde f \,
dx \in (- \infty, + \infty)
 \ .
 \end{equation}

Concerning the second addendum, differentiating a composition of
functions shows immediately that
\begin{equation}\label{uno} \lim _{t \to
0 ^ +} \frac{J \big ( f_{\alpha (t)} \big ) - J (f)}{t} = \delta J
(f, g)\ .
\end{equation}

By combining (\ref{due}) and (\ref{uno}), it is straightforward to
conclude. Indeed, similarly as in the proof of Theorem \ref{diffe},
we may distinguish the two cases $J (\tilde f)
> J (f)$ and $J (\tilde f) = J (f)$.

If $J (\tilde f) > J (f)$, the limit in (\ref{due}) remains finite,
whereas the limit in (\ref{uno}) becomes $+ \infty$. Hence it holds
$$\lim _{t \to 0 ^ +} \frac{J \big ((1-t) \cdot f \oplus t \cdot g \big ) - J (f)   }{t} = \delta J (f, g)  = + \infty\ ,$$
and the thesis of the lemma holds true.

If $J (\tilde f) = J (f)$,  then $\tilde f = f$ $\H ^n$-a.e., so
that the r.h.s.\ of (\ref{due}) agrees with $-\delta J (f, f)$, and
the lemma follows summing up (\ref{due}) and (\ref{uno}).

It remains to get rid of the assumption $v(0)= 0$. In the general
case, we set as usual $$d:= v (0)\, , \qquad \tilde v (x) := v (x) -
d\, , \qquad  \tilde g(x) := e ^ {- \tilde v(x)}\ .$$ Since
$$(1-t) \cdot f \oplus t \cdot g = e ^ {-dt} \big ((1-t) \cdot f \oplus t \cdot \tilde g
\big )\ , $$ we have
$$\frac{J\big ((1-t) \cdot f \oplus t \cdot g  \big )-J(f)}{t}
=J(f)\,\frac{e^{-dt}-1}{t} + e^{-dt}\,\frac{J\big ((1-t) \cdot f
\oplus t \cdot \tilde g \big )-J(f)}{t}\, .
$$
By passing to the limit as $t \to 0 ^+$, since $\tilde v (0) \geq 0$
by construction, we obtain
$$\lim_{t \to 0 ^+}
\frac{J\big ((1-t) \cdot f \oplus t \cdot g  \big )-J(f)}{t}
 \geq -d J (f) + \delta J (f, \tilde g) - \delta J (f, f)\ . $$
To conclude, it is enough to observe that  $ -d J (f) + \delta J (f,
\tilde g) = \delta J (f, g)$ ({\it cf.} (\ref{2.3})). \qed

\bigskip
{\it Proof of Theorem \ref{teomink}}. By the Pr\'ekopa-Leindler
inequality, the function $\psi (t):= \log \big ( J ((1-t) \cdot f
\oplus t \cdot g)\big )$ is concave on $[ 0, 1]$ ({\it cf.}\ Remark
\ref{remPL}). In particular, it holds
\begin{equation}\label{prebasic}
\psi (t) \geq \psi (0) + t \big [ \psi (1) - \psi (0)]\qquad \forall
t \in [0,1]\ . \end{equation} As a consequence, the (right)
derivative of the function $\psi$ at $t= 0$ satisfies
\begin{equation}\label{basic}
\psi' (0) \geq \big [\psi (1) - \psi (0) \big]\ . \end{equation} By
Lemma \ref{jerison}, we have
$$\psi' (0) = \frac{\delta J (f, g) - \delta J (f, f)}{J (f)}\ .
$$
Therefore (\ref{basic}) can be rewritten as
$$\frac{\delta J (f, g) - \delta J (f, f)}{J (f)} \geq \log \Big ( \frac{J (g)}{J (f)} \Big
)\ .$$ Inserting (\ref{forment}) into the above inequality,
(\ref{mink1}) is proved.

Finally, assume that $g (x) = f ( x-x_0)$ for some $x_0 \in \R ^n$.
Then (\ref{mink1}) holds with equality sign thanks to Proposition
\ref{ff} and the invariance of $J$ by translation of coordinates.
Conversely, assume that (\ref{mink1}) holds with equality sign. By
inspection of the above proof one sees immediately that also
inequality (\ref{basic}), and hence inequality (\ref{prebasic}),
must hold with equality sign. This entails that the
Pr\'ekopa-Leindler inequality holds as an equality, and therefore
$f$ and $g$ agree up to a translation. \qed

\section{Isoperimetric and log-Sobolev inequalities for log-concave
functions}\label{seclogsob} Let us now turn attention to some
consequences of the results in Sections \ref{sec2} and
\ref{secmink}.

Motivated by the equality
$$\lim_{t\to 0 ^+} \frac{V (K+tB_1)- V (K)} {t} =  P (K)\ ,$$
and having in mind that the Gaussian probability density
$$\g (x):= c_n e ^ { -  \frac{\| x\| ^ 2 }{2}}\ , \qquad c_n:=(2 \pi) ^ {-
\frac{n}{2}}\ ,$$ plays within the class $\A'$ the role of the unit
ball in $\K ^n$, we set the following

\begin{definition}\label{defper} {\rm For any $f \in \A'$ with $J (f) >0$,
we define the {\it perimeter} of $f$ as
$$P (f) :=  \delta J (f, \g)\ . $$ }
\end{definition}

Similarly as Minkowski first inequality (\ref{M1}) (when applied
with $L$ equal to a ball $B$) implies the classical isoperimetric
inequality
$${V(K) ^ {\frac{1}{n}}  } {P(K) ^ {-\frac{1}{ n-1}}} \leq {V(B) ^
{\frac{1}{n}}  } {P (B)^ {- \frac{1}{ n-1}}} \qquad \forall K \in
\K_0 ^n \ ,
$$
Theorem \ref{teomink} (when applied with $g = \g$ and combined with
Theorem \ref{A'}) yields the following functional version of the
isoperimetric inequality:

\begin{prop}\label{repper} Let $f  = e ^ {-u} \in \A'$,  and assume that $\varphi:= u^*$ is uniformly strictly convex, namely
\begin{equation}\label{HG}
\exists \, c>0 \ :\ \nabla ^2\varphi(y) \geq  c\, {\rm Id} \quad
\forall y \in \R ^n\ .
\end{equation}
Then
\begin{equation}\label{intrepper}
P (f) = \frac{1}{2} \int_{\sR ^n} \frac{ \|\nabla f\| ^ 2} { f} \, dx + (\log
c_n ) \, J (f) \geq n J (f) + \Ent (f) \ ,
\end{equation}
with equality sign if and only if there exists $x_0 \in \R ^n$ such
that $f (x) = \g  (x- x_0)$ $\forall x \in \R ^n$.
\end{prop}
\proof By the assumption (\ref{HG}), $\g$ is an admissible
perturbation for $f$ according to Definition \ref{defH}. Then, by
Definition \ref{defper} and Theorem \ref{A'}, one gets
$$P (f) = \delta J (f, \g) = \int_{\sR^n} ( \frac{1}{2} \| \nabla u \| ^ 2 + \log c_n )   f \, dx\ ,$$
which proves the first equality in (\ref{intrepper}). The subsequent
inequality in (\ref{intrepper}) is obtained by applying Theorem
\ref{teomink} (simply take into account that $J (\g) = 1$). \qed

\bigskip
As a further application of our results, we now provide a
generalized logarithmic Sobolev inequality for log-concave measures.
After the pioneering result by Gross concerning the Gaussian measure
\cite{Gross}, the validity of logarithmic Sobolev inequalities for
more general probability measures, having in particular a
log-concave density, has been investigated by several authors. We
refer in particular to the paper \cite{Bobkov} by Bobkov, where
necessary and sufficient conditions are discussed.

\begin{prop}\label{logsobolev}
Let $\nu = g \H ^n= e ^ {-v} \H ^n$ be a log-concave probability
measure such that $g \in \A'$ and
\begin{equation}\label{HG2}
 \nabla ^2 v  \geq  c\, {\rm Id} \quad
\hbox{ for some } c>0\ .
\end{equation}

Let $a: \R_+ \to \R_+$ be a continuous increasing function with
$a(0) = 0$.

Let $h$ be a positive measurable function of class $\C^2(\R ^n)$
which satisfies the conditions
$$\lim_{\|x \| \to + \infty} \frac{- \log (a(h))+v}{\|x\|} = + \infty
\quad \hbox{ and } \quad    - c' \, \nabla ^2 v  \leq \nabla ^ 2
(\log(a(h)) < \nabla  ^ 2 v \ \
 \hbox { for some } c'>0\ .
$$
Then it holds
\begin{equation}\label{LS}
 \int_{\sR ^n} a(h) \log a (h)
 d \nu   -  \Big ( \int_{\sR ^n} a(h) \, d \nu  \Big ) \log   \Big ( \int_{\sR ^n} a(h) \, d \nu  \Big )
 \leq  \frac{1}{c}
\int_{\sR ^n} \frac{(a' (h)) ^ 2} {a(h)} \|\nabla h \| ^ 2 \, d \nu
 \ .
 \end{equation}
\end{prop}

\begin{remark}{\rm The constant $\frac{1}{c}$ in the r.h.s.\ of (\ref{LS}) is non-optimal. Indeed, consider for instance the case when $g = \g$ (so that $c=1$), and $a (h ) = h ^2$. Then
 (\ref{LS}) becomes
\begin{equation}\label{LS'}
 \int_{\sR ^n} h ^2 \log  (h^2)
 d \nu   -  \Big ( \int_{\sR ^n} h^2 \, d \nu  \Big ) \log   \Big ( \int_{\sR ^n} h ^2 \, d \nu  \Big )
 \leq
4 \int_{\sR ^n} \|\nabla h \| ^ 2 \, d \nu
 \ ,
 \end{equation}
and it is known that (\ref{LS'}) holds true with $2$ in place of $4$ at the r.h.s. This assertion can be recovered by inspection of the proof below,
since in this case the number $t$ appearing in (\ref{st}) equals $\frac{1}{2}$.
 }
\end{remark}

\begin{remark}{\rm It is not surprising that, in order to have an inequality of logarithmic Sobolev type for the measure $\nu$, condition (\ref{HG2}) is needed; indeed, (\ref{HG2}) can be related to the so-called Herbst necessary condition (see \cite{Bobkov} for a more detailed discussion).}
\end{remark}

{\it Proof of Proposition \ref{logsobolev}.}
Since $\int _{\sR ^n} g \, dx = 1$, inequality (\ref{mink1})
reads
\begin{equation}\label{mink2}
\delta J (f, g) \geq n J (f) + \Ent (f) \ .
\end{equation}

The computation of $J (f)$ and $\int_{\sR^n} f \log f \, dx$ is
straightforward:
$$J (f)
=  \int_{\sR ^n }  a(h) \, d\nu  \ , \qquad  \int_{\sR^n} f \log f
\, dx =   \int_{\sR ^n } \big ( -  v + \log a (h) \big )  a(h) \,
d\nu\ . $$ On the other hand, by the hypotheses made on $h$ and $g$,
the functions $f:= a (h) g$ and $g$ turn out to satisfy the
assumptions of Theorem \ref{A'}. Then, setting $\psi = v ^*$, we
have
$$\delta J (f,
g) = \int _{\sR ^n} \psi  \big ( \nabla v - \nabla \log a (h) \big )
a (h) \, d \nu\ .
$$
Inserting the above expressions of $J (f)$, $\int_{\sR^n} f \log f
\, dx$, and $\delta J (f, g)$ into (\ref{mink2}) leads to
\begin{equation}\label{e1}
\int_{\sR ^n} a(h) \log a (h)
 d \nu   -  \Big ( \int_{\sR ^n} a(h) \, d \nu  \Big ) \log   \Big ( \int_{\sR ^n} a(h) \, d \nu  \Big )
 \leq R (h)\ ,
 \end{equation}
with $$ R(h) =  \int_{\sR^n}  \Big [  \psi \big (   \nabla v - \nabla
\log a (h)  \big )  + v - n \Big ] a (h) \, d \nu \ .$$ Using the
identity $v (x) = \langle x, \nabla v (x) \rangle - \psi (\nabla
v(x))$, we may rewrite $R (h)$ as
$$ R(h) =  \int_{\sR^n}  \Big [    \psi \big (\nabla v  - \nabla \log a (h)  \big ) - \psi (\nabla
v)  + \langle x, \nabla v  \rangle -n\Big ] a (h) \, d \nu \ .$$
Now we observe that
$$ \langle  x, \nabla v \rangle a (h ) g = - \langle  x, \nabla g\rangle a (h )
 = -{\rm div} ( x a(h) g) + \langle  x, \nabla a(h) \rangle g + {n}   a(h)g   \ ,$$ and
$$\int_{\sR^n} {\rm div} ( xa(h)g) \, dx = \lim _{ r \to + \infty} r \int_{\partial B _r}  a(h) g  \, d \H ^ {n-1}=
\lim _{ r \to + \infty} r \int_{\partial B _r} f \, d \H ^ {n-1} =
0$$ (where the last equality is satisfied by the exponential decay
of $f$ at infinity, {\it cf.}\ Lemma \ref{preli}). Therefore,
\begin{equation}\label{e2} R(h) = \int_{\sR^n}  \Big [     \psi \big (
 \nabla v - \nabla \log a (h) \big )- \psi (\nabla v)   + \langle  x, \nabla \log a(h)
\rangle \Big ] a (h) \, d \nu \ .\end{equation} In view of
(\ref{e1}) and (\ref{e2}), the statement is proved if the following pointwise
inequality holds:
$$ \psi \big ( \nabla v  - \nabla \log a (h) \big )- \psi
(\nabla v)   + \langle x, \nabla \log a(h)   \rangle  \leq
\frac{1}{c}
 \| \nabla \log a(h) \| ^ 2 \ .
$$
This is readily checked: indeed, setting $y := - \nabla \log a (h)$, by Lagrange theorem and
assumption (\ref{HG2}), there exist $t,s \in (0,1)$ such that
\begin{equation}\label{st}
\begin{array}{ll} \psi \big (  \nabla v +y  \big )- \psi
(\nabla v)   - \langle x , y \rangle & = \langle \nabla \psi (
\nabla v + ty), y \rangle - \langle \nabla \psi (\nabla v), y
\rangle   \\ \noalign{\smallskip} & =\displaystyle{ \langle \nabla ^
2 \psi ( \nabla v + sty) ty, y \rangle \leq \frac{1}{c}  \|y \| ^ 2}
 ,\end{array}
 \end{equation}
and the proof is achieved. \qed

\section{About the Minkowski problem}\label{towards}

In this concluding section we move the first steps towards the
solution of the functional Minkowski problem. In view of Theorems
\ref{A'} and \ref{A''}, its formulation within the class $\A'$ or
$\A''$ reads as follows:  find $f \in \A'$  such that
\begin{equation}\label{M'} \mu (f) = m\ ,
\end{equation} where $m$ is a given  positive Borel measure on $\R^n$, or find
$f \in \A''$ such that
\begin{equation}\label{M''}(\mu (f), \sigma (f)) = (m ,
\eta)\ ,\end{equation} where $(m, \eta)$ are given positive Borel
measures respectively  on $\R^n$ and $S ^ {n-1}$.
Here the measures $\mu (f)$ and $\sigma (f)$ are intended according
to Definition \ref{defmusigma}.

We begin by the following simple observation.

\begin{remark}\label{massafinita} {\rm We have the following finiteness necessary condition on
the measures $m$ and $\eta$, in order to solve the Minkowski problem with datum $m$ or $(m,\eta)$:
$$
\int_{\R^n} dm<+\infty\,,\quad
\int_{\sfe} d\eta<+\infty\,.
$$
Indeed, if $f$ belongs to $\A'$ or to $\A''$, we have
$$
\int_{\R^n}d\mu(f)=J(f)<+\infty\,,
$$
while, if $f\in\A''$ we have
$$
\int_{\sfe}d\sigma(f)\le (\max_{K} f) \,\H^{n-1}(\partial K)<+\infty\,,
$$
where $K=\dom(-\log f)$.
}
\end{remark}

\medskip

Next, we show that, for the solvability of (\ref{M'}), $m$
must satisfy an equilibrium condition, which is completely analogous
to the null barycenter property well-known in the classical
Minkowski problem for convex bodies. The same holds true, for the
solvability of (\ref{M''}), replacing $m$ by the pair $(m, \eta)$.

\begin{prop}\label{neccond}{

(i) For any $f\in \A'$, the measure $\mu (f)$ verifies
$$\int_{\sR^n} y\, d \mu (f) (y) = 0\ .$$

(ii) For any $f\in \A''$, the measures $\mu (f)$  and $\sigma (f)$
verify
$$\int_{\sR^n} y\, d \mu (f) (y) + \int_{S ^ {n-1}} y\, d \sigma (f) (y)  = 0\ .$$
}
\end{prop}

\proof Given a point $x_0 \in \R ^n$ and a function $v\in \L$, we
denote by $  [v]_{x_0}$ the translated function $x \mapsto v(x
-x_0)$. With this notation it is straightforward to check that, for
any $u, v \in \L$, it holds
\begin{equation}\label{transl}
u \Box  [v]_{x_o} =    [ u \Box v  ] _{x_0}\ .
\end{equation}

Assume now that $f = e ^ {-u}$ belongs either to $\A'$ or to $\A''$.
For any fixed $x_0 \in \R^n$ and any $\e >0$, let us compute $\delta
J (f, g_\e)$, where $g _\e = e ^ {-v_\e}$, being
$$v_\e (x):= \e u \big ( \frac{x-x_0}{\e} \big ) =
 [u \e] _{\frac{x_0}{\e}} (x) \qquad \forall \, x \in \R ^n\ .$$

For any $t>0$ one has
$$(v_\e t) =     [ u ( t \e)
 ]_{\frac{x_0}{\e}} \ , $$ and hence, in view of (\ref{transl}),
$$u \Box (v_\e t) = [ u \Box u ( t \e) ] _{\frac{x_0}{\e}}\ .$$
Therefore,
\begin{equation}\label{prima}
\delta J (f, g_\e) = \lim _{t \to 0 ^+} \frac{J ( e ^ {- u \Box u (
t \e)})  - J (f) }{t} = \e \lim _{t \to 0 ^+} \frac{J ( e ^ {- u
\Box u ( t \e)})  - J (f) }{ t\e } = \e \delta J (f, f)\ .
\end{equation}
On the other hand, we observe that
$$v_\e ^* (y) =  \langle x _0, y \rangle + \e u ^* (y)\qquad \hbox{ and } \qquad \dom (v_\e) =  x_0 + \e \dom (u) \ .$$
Therefore, if $f \in \A'$, by applying Theorem \ref{A'} we get
\begin{equation}\label{seconda'}
\delta J (f, g_\e) =   \int_{\sR ^n} \langle x _0, y \rangle  \, d
\mu (f) (y) +  \e \int_{\sR ^n} u ^* (y)  \, d \mu (f) (y) \ ;
\end{equation}
similarly, if $f \in \A''$, by applying Theorem \ref{A''}  we get
\begin{equation}\label{seconda''}
\begin{array}{ll} \delta J (f, g_\e) & \displaystyle{=   \int_{\sR ^n} \langle x _0, y \rangle
\, d \mu (f) (y) +  \e \int_{\sR ^n} u ^* (y)  \, d \mu (f) (y) }\\
\noalign{\smallskip} & \displaystyle{+ \int_{S ^{n-1}} \langle x _0,
y \rangle \, d \sigma (f) (y) + \e \int_{S ^{n-1}} h_{\dom (u)}  (y)
\, d \sigma (f) (y)}\ . \end{array}
\end{equation}
We now observe that the following terms, which appear multiplied by
$\e$ in (\ref{prima}), (\ref{seconda'}) and (\ref{seconda''}), are
finite:
$$\delta J (f, f)\ , \quad \int_{\sR ^n} u ^* (y)  \, d \mu (f)
(y)\ , \quad \int_{S ^{n-1}} h_{\dom (u)}  (y) \, d \sigma (f) (y)\
$$
(recall in particular Proposition \ref{ff} and Lemma
\ref{finitezza}). Then the statement follows by combining
(\ref{prima}) with (\ref{seconda'}) or (\ref{seconda''}), in the
limit as $\e \to 0 ^+$.
 \qed

\bigskip

\begin{remark}{\rm We observe that the conditions expressed by Remark \ref{massafinita} and Lemma
\ref{neccond} are in general not sufficient for the solvability of
the Minkowski problem within one of the classes $\A'$ or $\A''$.
Indeed, assume for instance that $n=1$ and consider the Minkowski
problem in $\A'$: given an absolutely continuous measure on $\R$
with a positive continuous density $m$, satisfying the necessary
conditions $\int _{\sR} m (y) \, dy < + \infty$ and $\int_{\sR} y m
(y) \, d y= 0$, it amounts to finding a function $\varphi \in \C ^ 2
_+(\R)$, with $u = \varphi ^* \in \L'$, solving the second order
o.d.e.
\begin{equation}\label{ode}
 e ^ {\varphi(y)  - y \varphi '(y) } \varphi ''(y) = m(y) \qquad \forall y \in
\R\ . \end{equation} We observe that, if $\varphi$ is a solution to
(\ref{ode}), for any $\alpha \in \R$, also $\varphi + \alpha y$ is a
solution. Therefore, we may assume with no loss of generality that
$\varphi' (0) = 0$, and write the unique solution to (\ref{ode})
with initial datum at $y=0$ as
\begin{equation}\label{sol}
\varphi (y) = \varphi (0) - y \int _0 ^y \frac { \log ( e ^
{\varphi(0)} - M (t)) - \varphi (0)}{t ^ 2 } \, dt \ , \qquad \hbox{
where } M (t) := \int _0 ^ t s m (s) \, ds\ .
\end{equation}
Now, in order that $u = \varphi ^* \in \L'$, we have to impose that
$\frac{\varphi (y)}{y}$ diverges as $|y| \to + \infty$. Such
condition can be satisfied (by inspection of (\ref{sol})) only if
\begin{equation}\label{value}
e ^ {\varphi (0)} = M _\infty:= \int_ 0 ^ {+ \infty} s m (s) \, ds\
. \end{equation} By (\ref{sol}) and (\ref{value}), it holds
$$\lim _{y \to + \infty} \frac{\varphi (y)}{y} =
\lim _{y \to + \infty} \int _0 ^y \frac { \log ( M _\infty - M (t))
- \log M _\infty }{t ^ 2 } \, dt\ .$$ It is quite easy to construct
explicit examples of
 positive continuous functions $m$, with finite integral and zero barycenter, such that limit at the r.h.s.\ of the above equality
 remains
 finite. For such a datum $m$, the Minkowski problem does not
 admit solutions in $\A'$.}
\end{remark}

\bigskip
In view of the above Remark, and since in higher dimensions equality
(\ref{M'}) does not correspond any longer to an o.d.e., but rather
to a Monge-Amp\`ere type equation, proving a general existence
result for the functional Minkowski problem seems to be a quite
delicate task. On the other hand, as a consequence of Corollary
\ref{cormink}, we are able to prove that uniqueness (up to
translations) holds true, in both the cases of $\A'$ and $\A''$.

\bigskip
\begin{prop}\label{uni} Let $f_1, f_2 \in \A$ satisfy one of the following conditions:
\begin{equation}\label{c1}
f_i \in \A' \quad i= 1, 2, \qquad \hbox{ and } \qquad \mu (f_1) =
\mu (f_2)
\end{equation}
or
\begin{equation}\label{c2}
f_i \in \A'' \quad i= 1, 2, \qquad \hbox{ and } \qquad \mu (f_1) =
\mu (f_2) , \ \sigma (f_1 ) = \sigma (f_2)\ .
\end{equation}
Then there exists $x_0\in \R ^n$ such that $f_2(x) = f_1 (x-x_0)$.
\end{prop}
\proof Firstly notice that the equality $\mu (f_1) = \mu (f_2)$
implies $J (f_1) = J (f_2)$. Moreover the assumption $f_i \in \A'$
(or $f_i \in \A''$) implies that $J (f_i) >0$. If (\ref{c1}) holds,
by Theorem \ref{A'} one has
$$\delta J (f_1, g) = \delta J (f_2, g) \qquad \forall g \in \A''\ .$$
In particular,  taking $g = f_1$ or $g= f_2$, one sees that condition (\ref{incroci}) is satisfied.
Therefore, we are in a position to apply Corollary \ref{cormink}, and the statement follows.
If (\ref{c2}) holds, the proof is exactly the same by using Theorem \ref{A''} in place of Theorem \ref{A'}. \qed

\section{Appendix}

This appendix contains the proofs of  some results stated in Section \ref{sec2}, precisely all the preliminary lemmas used in the
 proof of Theorems \ref{A'} and \ref{A''}, and the claim made in Remark \ref{remdim1}.

\bigskip
{\it Proof of Lemma \ref{chiusura}.}
It is immediate to
check that the classes $\L'$ and $\L ''$ are closed by right
multiplication by a positive scalar. Let us show that each of them
is closed also by infimal convolution.

\smallskip (i) Let $u, v \in \L '$, set $\varphi := u ^*$, $\psi := v
^*$, and $w:= u \Box v$.

By Proposition \ref{primolemma} (iii), it holds $\dom (w)= \dom (u)
+ \dom (v) = \R^n$.

The condition of having a superlinear growth at infinity is equivalent
to the condition of being cofinite \cite[Proposition 3.5.4]{BoVa},
and the latter is clearly closed by infimal convolution in view of
the equality $w^* = \varphi + \psi$ holding by Proposition
\ref{primolemma} (iv). Therefore, $w$ has superlinear growth at
infinity.

Since $(\R ^n, u)$ and $(\R ^n, v)$ are convex functions of Legendre
type, with $u, v \in \C ^ 2 _+$, the mappings $\nabla u$ and $\nabla
v$ are $\C^1$ bijections from $\R^n$ to $\R^n$, with a nonsingular
Jacobian. Therefore also their inverse maps, which by Proposition
\ref{proplegendre} are precisely $\nabla \varphi$ and $\nabla \psi$,
are $\C^1$ bijections from $\R^n$ to $\R^n$, and the same holds true
for their sum. Hence $(\R^n, \varphi + \psi)$ is a convex function
of Legendre type, with $\varphi + \psi$ of class $\C ^ 2_+$. In
turn, this implies that the Legendre conjugate of $(\R^n, \varphi +
\psi)$, namely $(\R^n, w)$, is a convex function of Legendre type,
with $w$ of class $\C ^ 2_+$.

\smallskip

(ii) Let $u, v \in \L ''$, and set $K: = \dom (u)$, $L:= \dom (v)$,
$\varphi, \psi$, and $w$ as above.

By Proposition \ref{primolemma} (iii), it holds $\dom (w)= K+ L \in
\K ^n \cap \C ^ 2_+$.

Since $u$ and $v$ are of class $\C ^ 2 _+$, and their gradients
diverge at the boundary of their domains,
 $({\rm int} (K), u)$ and $({\rm int} (L), v)$ are convex functions of Legendre
 type, and
the mappings $\nabla u$ and $\nabla v$ are $\C^1$ bijections
respectively from $K$ and $L$ onto $\R^n$. Hence, similarly as
above, we may apply Proposition \ref{proplegendre} to infer that
$(\R ^n, \varphi)$, $(\R ^n , \psi)$, and hence $(\R ^n, \varphi +
\psi)$, are convex functions of Legendre type, with $\varphi + \psi$
of class $\C ^ 2_+$. This yields that $(\R^n, w)$ is a convex
function of Legendre type, with $w$ of class $\C ^ 2_+$.

It remains to check that $w$ is continuous up to $\partial (K +L)$.
To this end we are going to use as a  crucial tool the identity
\begin{equation}\label{gen}
u \Box v (x) = \inf_{ x_1+x_2= x} \!\{ u (x_1) + v
(x_2)\} = u \big (\nu _K ^ {-1} (\nu _{K +L} (x)\big ) + v \big (\nu
_K ^ {-1} (\nu _{K +L} (x)\big ) \quad \forall x \in {\partial} ( K
+L)\ ,
\end{equation}
which follows from the definition of infimal convolution and the
assumption $\partial K,
\partial L \in \C ^ 2 _+$.

Let $\ov x\in
\partial (K +L)$, and let us show that for every sequence of points  $x ^ h \in K
+L$ such that $x ^h \to \ov x$, it holds
\begin{equation}\label{convbordo}
\lim_h u \Box v ( x ^ h) =   u \Box v  (\ov x) \ .
\end{equation}

Up to passing to a (not relabeled) subsequence, we may assume that
one of the following two cases occurs: $$x ^ h \in
\partial (K + L )
 \ \forall h \qquad \hbox{  or  } \qquad x ^h \in {\rm int} (K +L) \
 \forall h\ .$$

Consider first the case $x ^h \in \partial (K +L)$ $\forall h$. Let
us write the identity (\ref{gen}) at $ x ^h$
$$u \Box v (x^h) = u \big (\nu _K ^ {-1} (\nu _{K +L} (x^h)\big ) + v \big (\nu
_K ^ {-1} (\nu _{K +L} (x^h)\big ) \quad \forall h\ , $$ and then
let us pass to the limit in $h$. Since
 by hypothesis the Gauss maps $\nu _K$, $\nu _L$  and their inverse are continuous,
and $u$, $v$ are continuous  up to $\partial K$, $\partial L$, we
get
$$\lim _h u \Box v (x^h) =u \big (\nu _K ^ {-1} (\nu _{K +L} (\ov x)\big ) + v \big (\nu
_K ^ {-1} (\nu _{K +L} (\ov x)\big )\ .$$ In view of the identity
(\ref{gen}), the r.h.s. of the above equality equals $u \Box v (\ov
x)$, and (\ref{convbordo}) is proved.

Consider now the case $x ^h \in {\rm int} (K +L)$ $\forall h$. We
set
$$
y ^h := \nabla w (x^h)= (\nabla (\varphi + \psi))^ {-1} (x^h)\ ,$$
and we decompose $x ^h$ as $x_1^h + x _2 ^h$, with
$$ x^h _1 := \nabla \varphi (y ^h)\in {\rm int} (K) \qquad
\hbox{and} \qquad x^h _2 := \nabla \psi (y ^h) \in {\rm int}(L)\ .$$
Then we have
$$u \Box v ( x ^ h) = [\langle x _1 ^ h , y ^ h \rangle - \varphi (y ^h)]  + [\langle x _2 ^ h , y ^ h \rangle - \psi (y ^h)]
 = u ( x _ 1 ^h)  + v ( x_ 2 ^h)\ .
$$
Let us now pass the the limit in $h$. By compactness, after possibly
selecting a (not relabeled) subsequence, there exist $\lim _h x _ 1
^h =: \ov x_1 \in
\partial K$ and $\lim _h x _ 2 ^h =: \ov x_2 \in
\partial L$. Since by assumption $u \in \C ^0 (K)$ and $v \in \C ^0 (L)$, we infer
$$\lim_h u \Box v ( x^ h) = u ( \ov x _ 1)  + v ( \ov x_ 2 )\
$$
In view of the identity (\ref{gen}), the above equality implies
(\ref{convbordo}) provided
$$\ov x_1 = \nu _K ^ {-1} \big (\nu _{K +L} (\ov x) \big ) \qquad \hbox  { and } \qquad
\ov x_2 = \nu _L ^ {-1} \big (\nu _{K +L} (\ov x) \big )\ .$$ In
turn, by the $\C ^ 2 _+$ assumption on $\partial K, \partial L$,
such conditions are satisfied provided the normal vectors $\nu _K
(\ov x_1)$ and $\nu _L (\ov x_2)$ coincide.
 Let us show that in fact each of them agrees with
$$\ov \xi := \lim _h \frac{y ^ h} {\| y ^ h \|}\ .$$
Since $y ^ h = \nabla u ( x  ^ h _1)$, and $\| y ^ h \| \to +
\infty$ (being $y ^h = \nabla w_h (x^h)$ and $x ^ h \to \ov x \in
\partial (K +L)$), by passing to the limit in the inequality
$$\frac{u (x) }{\|y ^h\|}  \geq  \frac{u (x^h _1) }{\|y ^h\|}  + \langle \frac{y ^ h} {\| y ^ h \|} , x - x ^ h _1 \rangle\ ,$$
we infer that any cluster point of the sequence $y ^ h  / \| y ^h\|$
belongs to the normal cone to $\partial K$ at $\ov x _1$, which is
reduced to $\nu _K (\ov x_1)$. In the same way we obtain $\ov \xi =
\nu _L (\ov x_2)$, and the proof is achieved. \qed

\bigskip\medskip
{\it Proof of Lemma \ref{convergenza}.}
(i) Let $x \in  \dom (u)$ be fixed. By the assumption $v(0)= 0$, we
have $u _t(x) \leq u(x)$ for every $t >0$, so that $\limsup _{t
\to 0^+} u _t (x) \leq u (x)$. Let us prove that we also have
\begin{equation}\label{limi}
\liminf_{t \to 0^+} u _t (x) \geq u (x)
\end{equation}
Assume $u, v \in \L'$, and set $\varphi := u ^*$, $\psi := v ^*$.
We choose $r >
\|\nabla u (x)\|  $ and we set $c:= \sup _{B _r} \psi$ (notice that
$c$ is finite because $\psi$ is
bounded on bounded sets \cite[Theorem 4.4.13]{BoVa}).  Then
$$\begin{array}{ll} u _t (x) & \displaystyle{= \sup _{y \in \sR ^n } \big \{ \langle x, y \rangle - \varphi (y)
- t \psi (y) \big \} \geq \sup _{y \in B _r } \big \{ \langle x, y
\rangle - \varphi (y) \big \}  -tc }  \\ \noalign{\smallskip} & =
\langle x, \nabla u (x) \rangle - \varphi (\nabla u (x)) -tc = u (x)
-  t c \, ,
\end{array}
$$
and (\ref{limi}) follows by passing to the inferior limit as $t \to 0 ^+$.

Assume $u, v \in \L''$. Setting $L:= \dom (v)$ and $m:= \min v$, it holds $v \geq I _L
+m $. Then
$$\begin{array}{ll} u _t (x) & \displaystyle{= \inf  _{x_1 + x _2 = x } \big \{u (x_1) + t v (x_2/t) \big \}  \geq \inf _{x_1 + x _2 = x }
\big \{ u (x_1) + tI _L (x_2/t) \big \} + t m  }
\\ \noalign{\smallskip} & =
\displaystyle{\inf _{x_1 + x _2 = x }\big \{ u (x_1) + tI
_{tL}(x_2)\big \} + t m =  \inf _{x_1 \in K \cap (x-  tL)  } \{ u
(x_1) \} + t m}
 \, ,
\end{array}
$$
and, thanks to the continuity of $u$ at $x$, (\ref{limi}) follows by passing to the inferior limit as $t \to 0 ^+$.

\smallskip
Statement (ii) is an immediate consequence of the convexity of the
functions $u _t$ and of the differentiability of their pointwise
limit $u$ in the interior of its domain. \qed

\bigskip
{\it Proof of Lemma \ref{puntuale}.}
Set $K _t := \dom (u _t)$.
First we claim that, for every fixed $x \in {\rm int} (K _t)$,
\begin{equation}\label{differe}
\hbox{ the map }\  t \mapsto \nabla u _t (x) \ \hbox{ is differentiable on $(0, + \infty)$. }
\end{equation}
Indeed, as noticed in the proof of Lemma \ref{chiusura}, the Fenchel conjugates $\varphi:= u ^*$ and $\psi:= v ^*$
are both of class $\C ^2_+$ on $\R^n$. Therefore,  the function
$F: \R ^n \times \R ^n \times (0, + \infty) \to \R ^n$ defined by
$$F (x, y, t):= \nabla \varphi (y) + t \nabla \psi (y) - x\ ,$$
is of class $\C ^1$ on $\R ^n \times \R ^n \times (0, + \infty)$, and
$\frac{\partial F}{\partial y} = {\nabla ^ 2 \varphi} + t \nabla ^ 2 \psi$ is nonsingular for every $y \in \R ^n$. Consequently, by the implicit function theorem, the equation $F (x, y, t)= 0$ locally defines a map $y = y (x, t)$ which is of class $\C ^1$ in its arguments.
By Lemma \ref{chiusura},
$( {\rm int} (K_t), {u_t})$ is a convex function of Legendre type, hence by Proposition \ref{proplegendre} $\nabla u _t$ is the inverse map of $\nabla \varphi _t$, namely
$$F (x, \nabla u _t (x), t) = \nabla \varphi _t (\nabla u _t (x)) -x = 0\ .$$
Therefore, for every $x \in {\rm int} (K_t)$ and every $t>0$, $y (x, t) = \nabla u _t (x)$, and (\ref{differe}) is proved.

Next, we apply again to Proposition \ref{proplegendre} in order to write the identity
\begin{equation}\label{identita}
u _t (x) =\langle x , \nabla u _t(x) \rangle -
\varphi _t \big ( \nabla u _t (x) \big ) \qquad \forall x \in {\rm int} (K_t)\ .
\end{equation}

By (\ref{differe}) and (\ref{identita}) we obtain that, for every fixed $x \in {\rm int} (K_t)$, the map $t \mapsto u _t (x)$ is differentiable on $(0, + \infty)$, with
$$\frac{d}{dt} u _t (x) = \langle x ,  \frac{d}{dt} \big ( \nabla u _t (x) \big )\rangle - \psi   \big ( \nabla u _t (x) \big )
- \langle \nabla \varphi _t  \big ( \nabla u _t (x) \big ) , \frac{d}{dt}  \big ( \nabla u _t (x) \big )\rangle  =
- \psi  \big ( \nabla u _t (x) \big )\ .$$

\qed

\bigskip
{\it Proof of Lemma \ref{finitezza}.}
We have
$$\begin{array}{ll}\displaystyle{\int_{\sR^n} \varphi (\nabla u(x) ) \, f(x) \, dx }& = \displaystyle{\int_{\sR ^n} \big ( \langle x
, \nabla u  \rangle - u   \big )\, f  \, dx  = - \int _{\sR ^n}
\langle x , \nabla f  \rangle \, dx +\int_{\sR^n} f \log f \, dx\ .}
\\
\noalign{\medskip}  & = \displaystyle{ - \int _{\sR ^n} {\rm div} (
f x ) \, dx + n J (f) + \int_{\sR^n} f \log f \, dx\ . }
\end{array}
$$
We observe that
$$\int _{\sR
^n} {\rm div} ( f x ) \, dx  = \lim _{ r \to + \infty} \int _{B_r }
{\rm div} ( f x ) \, dx =  \lim _{ r \to + \infty} r \, \int
_{\partial B_  r }  f  \, d {\mathcal H} ^ {n-1}= 0  \ ,$$ where the
last equality holds true by Lemma \ref{preli}. Therefore we have
$$\int_{\sR^n} \varphi (\nabla u(x) ) \, f(x) \, dx  = n J (f) +\int_{\sR^n} f \log f \, dx\ ,$$
and the lemma follows recalling that both $J (f)$ and $\int_{\sR^n}
f \log f \, dx$ are finite ({\it cf.}\ respectively Lemma
\ref{preli} and Proposition \ref{ff}).
 \qed

\bigskip
{\it Proof of Lemma \ref{ultimo}.}
By definition we have
$$
u _t (x)  = \inf _{x _1 + x _2 = x} \Big \{ u ( x_1) + t v \big ( \frac{x_2}{t} \big ) \Big \}\ .$$
Since
$$v_{\min} + I _L (x) \leq v(x) \leq v_{\max} + I _L (x) \qquad \forall x \in \R^n\ ,$$
it holds
$$\inf _{x _1 + x _2 = x} \Big \{ u ( x_1) + t v _{\min}  + t I _L\big ( \frac{x_2}{t} \big )\Big \}
\leq  u _ t (x)  \leq
\inf _{x _1 + x _2 = x} \Big \{ u ( x_1) + t v _{\max}  + t I _L\big ( \frac{x_2}{t} \big )\Big \}\ ,
$$
namely
$$ t v _{\min}  + \inf _{x _1 \in K \cap ( x- tL) } \big \{ u ( x_1) \big \}
\leq  u _ t (x)  \leq t v _{\max}  +
\inf _{x _1 \in K \cap ( x- tL) } \big \{ u ( x_1) \big \}\ .
$$
Therefore the statement is satisfied by taking $y$ as a point where $u$  attains its minimum on $K \cap ( x- tL )$. \qed

\bigskip
{\it Proof of Remark \ref{remdim1}.}
By inspection of the proof of Theorems \ref{A'} and \ref{A''}, one can see that assumption (\ref{H}) is used only in Step 4 (in order to prove that $\Psi (0) < + \infty$)
and in Step 5 (in order to prove that $\lim _{s \to 0^+} \Psi (s) = \Psi (0)$). Assume now $n=1$, and drop assumption (\ref{H}): let us indicate how Steps 4 and 5 (and consequently also Step 6) have to be modified in order to show that (\ref{tesi''}) continues to hold, possibly as an equality $+ \infty = + \infty$.

In Step 4, we limit ourselves to prove that $\Psi$ takes finite values at every $s>0$.

In Step 5, the proof of the continuity of $\Psi$ at every $s>0$ remains unchanged, whereas for $s\to 0 ^+$ we make the following claim (whose proof is postponed below):
\begin{equation}\label{onedim}
\hbox{ if $\Psi (0) < + \infty$, then $\Psi$ is continuous from the right at $s=0$}.
\end{equation}

Consequently, in Step 6 we must distinguish two cases. In case $\Psi (0) < + \infty$, thanks to (\ref{onedim}) equality (\ref{interno}) can be proved exactly as before. In case $\Psi (0) = + \infty$, (\ref{interno}) continues to hold as an equality $+ \infty = + \infty$, and it can be proved by slight modifications of the case $\Psi (0) <  \infty$. More precisely, (\ref{positives}) and (\ref{limder}) in Step 6 remain unchanged, whereas (\ref{caso1})
has to be replaced by
\begin{equation}\label{caso2}
\liminf _{s \to 0 ^+} \Psi (s) \geq
\sup _{C \subset \subset E} \liminf _{s \to 0 ^+} \int _{C} \psi \,
d \mu (f_s)  = \sup _{C \subset \subset E} \int _{C} \psi  d \mu (f)
= + \infty \end{equation}
(notice that the second equality in (\ref{caso2}) holds by dominated convergence, since  by Lemma \ref{convergenza} we have
$\psi (\nabla u _s) f _s \to \psi (\nabla u) f$ as $s \to 0 ^+$,  and   by
Lemma \ref{mono} and Lemma \ref{convergenza} (ii)
the nonnegative functions $ \psi (\nabla u _{s})  \,
f_{s}$ are bounded above on $C$ by some continuous function
independent of $s$).

Let us finally prove (\ref{onedim}). Assume
\begin{equation}\label{psi0}
\Psi (0) = \int _{\sR} \psi \varphi '' e ^ {\varphi - y \varphi'} \, dy < + \infty\ .
\end{equation}
Since $n=1$, (\ref{scompo}) simplifies into
$$\Psi (s) = \Psi _0 (s) + s \Psi _1 (s)\ ,$$
where
$$\begin{array}{ll}
& \displaystyle {\Psi _0(s): = \int_{\sR} h _0 (s, y) \, dy \, \qquad h _0 (s, y):=
\psi e ^ {\varphi _s  - y \varphi' _s
 } \varphi'' \chi _{Q_s}}  \\
& \displaystyle {\Psi _1(s): = \int_{\sR} h _1 (s, y) \, dy \, \qquad h _1 (s, y):=
\psi  e ^ {\varphi _s  - y \varphi' _s
 } \psi'' \chi _{Q_s}\ .}
\end{array}$$
To get (\ref{onedim}) it suffices to show that
\begin{eqnarray}
& \lim _{s \to 0 ^ +} \Psi _0 (s) = \Psi (0)  & \label{I'}
\\ & \lim_{s \to 0 ^ +} s \Psi _1 (s) =0 \ .  & \label{II'}
\end{eqnarray}
Thanks to assumption (\ref{psi0}), (\ref{I'}) can be proved exactly as before ({\it cf.} the proof of (\ref{I})).
To prove (\ref{II'}), we write
$$s \Psi _1 (s) = I _+ (s) + I _- (s) := \int _{\sR _+} s h _1 (s, y) \, dy +
\int _{\sR _-} s h _1 (s, y) \, dy \ ,$$
and we show that both $I _\pm (s)$ are infinitesimal as $s \to 0^+$. Let us consider $I _+ (s)$ (the case of $I _- (s)$ is completely analogous).

We observe that
$$0 \leq sh _1 (s, y) = s \psi  e ^ {\varphi _s  - y \varphi' _s
 } \psi'' \chi _{Q_s} \leq - F (y) G' _s (y) \qquad \forall y, s>0\ ,$$
where we have set
$$F(y) := \frac{ \psi}{y} e ^ {\varphi - y \varphi '} \qquad \hbox{ and } \qquad G_s (y) := e ^ {s (\psi - y \psi')}\ .$$
Then an integration by parts gives
$$0 \leq I _ + (s) \leq \lim _{\varepsilon \to 0^+, \ r \to + \infty} \Big \{\int_{\varepsilon} ^ r F' (y) G_s (y) \, dy + F (\varepsilon) G_s (\varepsilon)- F (r) G_s (r)  \Big \}\ .$$
Since $\psi(0) = \psi ' (0) = 0$ (respectively because $v\geq  0$ and $\psi \geq 0$), passing to the limit in $\varepsilon$ gives
\begin{equation}\label{provvisoria}
0 \leq I _ + (s) \leq \lim _{\ r \to + \infty} \Big \{\int_{0} ^ r F' (y) G_s (y) \, dy - F (r) G_s (r)  \Big \}\ .
\end{equation}
Next we observe that the following limit exists:
$$\alpha := \lim_ {r \to + \infty} F(r) = \lim_ {r \to + \infty} \int_0 ^ r F' (y) dy\ .$$
Indeed a straightforward computation gives
$$F' (y) = - \psi \varphi '' e ^ {\varphi - y \varphi '} + \frac{ e ^ {\varphi - y \varphi '} }{y ^2} (\psi' y - \psi)\ ,$$
and both the functions at the right and side are integrable on $(0, + \infty)$ (the former by assumption (\ref{psi0}), the latter because it is nonnegative).

Let us show that $\alpha>0$ cannot occur. Indeed in such case, for some constants $\ov c$ and $\ov r$, it would be $F(r) \geq \ov c$  $\forall r \geq \ov r$. This would contradict (\ref{psi0}), since
$$\Psi (0) \geq \int _{\ov r} ^ {+ \infty} \psi \varphi '' e ^ {\varphi - y \varphi'} \, dy
\geq \ov c \int_{ \ov r} ^ { + \infty} y \varphi '' \, dy = \ov c \, \Big \{ \lim _{r \to + \infty} [ r \varphi ' (r) - \varphi (r) ]  -
 [ \ov r \varphi ' (\ov r) - \varphi (\ov r) ] \Big \}  = + \infty\ .$$

Taking into account  that $\alpha = 0$  (and also that $\lim _{r \to + \infty } G_s (r) = 0$), we may rewrite (\ref{provvisoria}) as
\begin{equation}\label{riprovvisoria}
0 \leq I _ + (s) \leq \lim _{\ r \to + \infty} \int_{0} ^ r F' (y) G_s (y) \, dy \ .
\end{equation}
Moreover, since $\alpha = 0$, we have in particular $\int _0 ^ {+ \infty} F'(y)\, dy < + \infty$, which implies  $F' \in L ^ 1 (0, + \infty)$.
Therefore, for every fixed $s>0$, the functions $F' G_s$  satisfy
$$|F' (y) G_s (y) | \leq |F'(y)| \in L ^ 1 (0, + \infty)\ .$$
We deduce that (\ref{riprovvisoria})  can be rewritten as
\begin{equation}\label{definitiva}
0 \leq I _+ (s) \leq \int _0 ^ {+ \infty} F' (y) G_s (y) \, dy\ .
\end{equation}
Finally, passing to the limit as $s\to  0 ^+$ in the right  hand side of (\ref{definitiva}) we obtain
$$\lim _{s \to 0 ^ +} \int _0 ^ {+ \infty} F' (y) G_s (y)  \, dy =
 \int _0 ^ {+ \infty} \lim _{s \to 0 ^ +} F' (y) G_s (y)  \, dy =  \int _0 ^ {+ \infty} F' (y)   \, dy  = \alpha = 0\ .$$
This implies that $I _+ (s)$ is infinitesimal as $s \to 0 ^+$ and the proof is achieved.
\qed

\end{document}